\documentclass[12pt]{article}
\usepackage{amsmath}
\usepackage{amsfonts}
\usepackage{amssymb}
\usepackage{color}
\usepackage{enumerate}

\newtheorem{thm}{Theorem}[section]
\newtheorem{cor}[thm]{Corollary}
\newtheorem{lem}[thm]{Lemma}

\newtheorem{rem}[thm]{Remark}
\providecommand{\U}[1]{\protect\rule{.1in}{.1in}}
\begin{document}

\title{Impulse output rapid stabilization for heat equations}
\author{Kim Dang Phung\thanks{
Universit\'{e} d'Orl\'{e}ans, Laboratoire MAPMO, CNRS UMR 7349, F\'{e}d\'{e}%
ration Denis Poisson, FR CNRS 2964, B\^{a}timent de Math\'{e}matiques, B.P.
6759, 45067 Orl\'{e}ans Cedex 2, France. E-mail address:
kim\_dang\_phung@yahoo.fr.} , Gengsheng Wang\thanks{
School of Mathematics and Statistics, and Collaborative Innovation Centre of
Mathematics, Wuhan University, Wuhan, 430072, China. E-mail address:
wanggs62@yeah.net. The author was partially supported by the National
Natural Science Foundation of China under grant 11571264.} , Yashan Xu%
\thanks{
School of Mathematical Sciences, Fudan University, KLMNS, Shanghai 200433,
China. E-mail address: yashanxu@fudan.edu.cn. }}
\date{}
\maketitle

\begin{abstract}
The main aim of this paper is to provide a new feedback law for the heat
equations in a bounded domain $\Omega $ with Dirichlet boundary condition.
Two constraints will be compulsory: First, The controls are active in a
subdomain of $\Omega $ and at discrete time points; Second, The observations
are made in another subdomain and at different discrete time points. Our
strategy consists in linking an observation estimate at one time, minimal
norm impulse control, approximate inverse source problem and rapid output
stabilization.
\end{abstract}

\vspace{0.2cm}

 \vspace{0.2cm}


\bigskip

\bigskip

\section{Introduction}

\bigskip

\bigskip

Let $\Omega \subset \mathbb{R}^{d}$, $d\geq 1$, be a bounded domain, with a $%
C^{2}$ boundary $\partial \Omega $. Let $V$ be a function in $L^{\infty
}\left( \Omega \right) $ with its norm $\left\Vert \cdot \right\Vert
_{\infty }$. Define 
\begin{equation*}
A:=\Delta -V\text{, with }D\left( A\right) =H^{2}\left( \Omega \right) \cap
H_{0}^{1}\left( \Omega \right) \text{ .}
\end{equation*}%
Write $\left\{ e^{tA},t\geq 0\right\} $ for the semigroup generated by $A$
on $L^{2}\left( \Omega \right) $. It is well-known that when $V=0$, the
semigroup $\left\{ e^{tA},t\geq 0\right\} $ has the exponential decay with
the rate $\alpha _{1}$, which is the first eigenvalue of $-\Delta $ with the
homogeneous Dirichlet boundary condition. The aim of this study is to build
up, for each $\gamma >0$, an output feedback law $\mathcal{F}_{\gamma }$ so
that any solution to the closed-loop controlled heat equation, associated
with $A$ and $\mathcal{F}_{\gamma }$, has an exponential decay with the rate 
$\gamma $. Further, two constraints are imposed: We only have accesses to
the system at time $\frac{1}{2}T+\overline{\mathbb{N}}T$ on an open
subdomain $\omega _{1}\subset \Omega $; We only can control at time $T+%
\overline{\mathbb{N}}T$ on another open subset $\omega _{2}\subset \Omega $.
Here, $T$ is an arbitrarily fixed positive number and $\overline{\mathbb{N}}%
:=\mathbb{N}\cup \left\{ 0\right\} $, with $\mathbb{N}$ the set of all
positive natural numbers. Therefore, the closed-loop controlled equation
under consideration reads:%
\begin{equation}
\left\{ 
\begin{array}{ll}
y^{\prime }\left( t\right) -Ay\left( t\right) =0\text{ ,} & t\in \left(
0,\infty \right) \backslash \mathbb{N}T\ \text{,} \\ 
y\left( 0\right) \in L^{2}\left( \Omega \right) \text{ ,} &  \\ 
y\left( nT\right) =y\left( nT_{-}\right) +1_{\omega _{2}}\mathcal{F}_{\gamma
}\left( 1_{\omega _{1}}^{\ast }y\left( \left( n-\frac{1}{2}\right) T\right)
\right) \text{ ,} & \forall n\in \mathbb{N}\text{ .}%
\end{array}%
\right.  \tag{1.1}  \label{1.1}
\end{equation}%
Here, $y\left( nT_{-}\right) $ denotes the left limit of the function: $%
t\mapsto y\left( t\right) $ (from $\left[ 0,\infty \right) $ to $L^{2}\left(
\Omega \right) $) at time $nT$; $1_{\omega _{2}}$ denotes the zero-extension
operator from $L^{2}\left( \omega _{2}\right) $ to $L^{2}\left( \Omega
\right) $ (i.e., for each $f\in L^{2}\left( \omega _{2}\right) $, $1_{\omega
_{2}}\left( f\right) $ is defined to be the zero-extension of $f$ over $%
\Omega $); $1_{\omega _{1}}^{\ast }$ stands for the adjoint operator of $%
1_{\omega _{1}}$; $\mathcal{F}_{\gamma }$\ is a linear and bounded operator
from $L^{2}\left( \omega _{1}\right) $ to $L^{2}\left( \omega _{2}\right) $.
The operator $\mathcal{F}_{\gamma }$ is what we will build up. The evolution
distributed system (\ref{1.1}) is well-posed and can be understood as the
coupling of a sequence of heat equations: 
\begin{equation*}
y\left( t\right) :=\left\vert 
\begin{array}{ll}
y^{0}\left( t\right) \text{ ,} & \text{if }t\in \left[ 0,T\right) \\ 
y^{n}\left( t\right) \text{ ,} & \text{if }t\in \left[ nT,\left( n+1\right)
T\right)%
\end{array}%
\right.
\end{equation*}%
for any $n\in \mathbb{N}$, where 
\begin{equation*}
\left\{ 
\begin{array}{ll}
\begin{array}{ll}
\partial _{t}y^{0}-\Delta y^{0}+Vy^{0}=0\text{ ,} & \text{\ in }\Omega
\times \left( 0,T\right) \ \text{,} \\ 
y^{0}=0\text{ ,} & \text{\ on }\partial \Omega \times \left( 0,T\right) \ 
\text{,}%
\end{array}
&  \\ 
\text{~}y^{0}\left( 0\right) =y\left( 0\right) \in L^{2}\left( \Omega
\right) \text{ ,} & 
\end{array}%
\right.
\end{equation*}%
and%
\begin{equation*}
\left\{ 
\begin{array}{ll}
\begin{array}{ll}
\partial _{t}y^{n}-\Delta y^{n}+Vy^{n}=0\text{ ,} & \text{\ in }\Omega
\times \left( nT,\left( n+1\right) T\right) \ \text{,} \\ 
y^{n}=0\text{ ,} & \text{\ on }\partial \Omega \times \left( nT,\left(
n+1\right) T\right) \ \text{,}%
\end{array}
&  \\ 
\text{~}y^{n}\left( nT\right) =y^{n-1}\left( nT\right) +1_{\omega _{2}}%
\mathcal{F}_{\gamma }\left( 1_{\omega _{1}}^{\ast }y^{n-1}\left( \left( n-%
\frac{1}{2}\right) T\right) \right) \in L^{2}\left( \Omega \right) \text{ .}
& 
\end{array}%
\right.
\end{equation*}

\bigskip

Throughout this paper, we denote by $\left\Vert \cdot \right\Vert $ and $%
\left\langle \cdot ,\cdot \right\rangle $ the norm and the inner product of $%
L^{2}\left( \Omega \right) $ respectively. The notation $\left\Vert \cdot
\right\Vert _{\omega _{1}}$ and $\left\langle \cdot ,\cdot \right\rangle
_{\omega _{1}}$ will mean the norm and the inner product of $L^{2}\left(
\omega _{1}\right) $ respectively. We write $\mathcal{L}\left(
H_{1},H_{2}\right) $ for the space consisting of all bounded linear
operators from one Hilbert space $H_{1}$ to another Hilbert space $H_{2}$.
Lastly, we introduce the set $\left\{ \lambda _{j}\right\} _{j=1}^{\infty }$
for the family of all eigenvalues of $-A$ so that%
\begin{equation}
\lambda _{1}\leq \lambda _{2}\leq \cdot \cdot \leq \lambda _{m}\leq
0<\lambda _{m+1}\leq \cdot \cdot \cdot \text{ and }\underset{j\rightarrow
\infty }{\text{lim}}\lambda _{j}=\infty \text{ ,}  \tag{1.2}  \label{1.2}
\end{equation}%
and let $\left\{ \xi _{j}\right\} _{j=1}^{\infty }$ be the family of the
corresponding orthogonal normalized eigenfunctions.

\bigskip

The main theorem of this paper will be precisely presented in section 5. It
can be simply stated as follows: For each $\gamma >0$, there is $\mathcal{F}%
_{\gamma }\in \mathcal{L}\left( L^{2}\left( \omega _{1}\right) ,L^{2}\left(
\omega _{2}\right) \right) $ and a positive constant $C_{\gamma }$
(depending on $\gamma $ but independent of $t$) so that each solution $y$ to
the equation (\ref{1.1}) satisfies the inequality:%
\begin{equation}
\left\Vert y\left( t\right) \right\Vert \leq C_{\gamma }e^{-\gamma
t}\left\Vert y\left( 0\right) \right\Vert \quad \text{\ for any }t\geq 0%
\text{ .}  \tag{1.3}  \label{1.3}
\end{equation}%
We now give two comments on this result. First, the aforementioned $\mathcal{%
F}_{\gamma }$ has the form:%
\begin{equation}
\mathcal{F}_{\gamma }\left( p\right) =-\sum_{j=1}^{K_{\gamma }}e^{\lambda
_{j}T/2}\langle g_{j},p\rangle _{\omega _{1}}f_{j}\left( x\right) \quad 
\text{\ for any }p\in L^{2}\left( \omega _{1}\right) \text{ .}  \tag{1.4}
\label{1.4}
\end{equation}%
Here, $K_{\gamma }\in \mathbb{N}$ is the number of all eigenvalues $\lambda
_{j}$ which are less that $\gamma +\frac{\text{ln}2}{T}$; $g_{j}$ and $f_{j}$
are vectors in $L^{2}\left( \omega _{1}\right) $ and $L^{2}\left( \omega
_{2}\right) $ respectively. These vectors are minimal norm controls for a
kind of minimal norm problems which can be given by constructive methods.
Second, the operator norm of $\mathcal{F}_{\gamma }$ is bounded by $%
C_{1}e^{C_{2}\gamma }$, with $C_{1}>0$ and $C_{2}>0$ independent of $\gamma $%
.

\bigskip

We next explain our strategy and key points to prove the above-mentioned
results. First, we realize that if a solution $y$ to the equation (\ref{1.1}%
) satisfies the inequality:%
\begin{equation*}
\left\Vert y\left( \left( \tfrac{5}{4}+n\right) T\right) \right\Vert \leq
e^{-\gamma T}\left\Vert y\left( \left( \tfrac{1}{4}+n\right) T\right)
\right\Vert \quad \text{\ for all }n\in \overline{\mathbb{N}}\text{ ,}
\end{equation*}%
then (\ref{1.3}) holds. Next, by the time translation invariance of the
equation, we can focus our study on the interval $\left[ \left( \frac{1}{4}%
+n\right) T,\left( \frac{5}{4}+n\right) T\right) $. And the problem of
stabilization is transferred into the following approximate controllability
problem (of the system (\ref{1.1}) over the above interval): Find a control
in a feedback form driving the system from each initial datum $y\left(
\left( \frac{1}{4}+n\right) T\right) $ to $y\left( \left( \frac{5}{4}%
+n\right) T\right) $ with the above estimate. When building up the feedback
law, we propose a new method to reconstruct approximatively the initial data 
$y\left( \left( \frac{1}{4}+n\right) T\right) $ from the knowledge of $%
1_{\omega _{1}}^{\ast }y\left( \left( \frac{1}{2}+n\right) T\right) $. All
along such process, we need to take care of the cost of $\mathcal{F}_{\gamma
}$. Naturally, we may expect that smaller is $e^{-\gamma T}$, larger is the
cost. It is worth mentioning, that by projecting the initial data $y\left(
\left( \frac{1}{4}+n\right) T\right) $ into subspaces span$\left\{ \xi
_{1},\cdot \cdot \cdot ,\xi _{K_{\gamma }}\right\} $ and span$\left\{ \xi
_{K_{\gamma }+1},\xi _{K_{\gamma }+2},\cdot \cdot \cdot \right\} $
respectively, we have 
\begin{equation*}
\left\Vert P_{K_{\gamma }}y\left( \left( n+1\right) T\right) \right\Vert
\leq e^{-\lambda _{K_{\gamma }+1}3T/4}\left\Vert y\left( \left( \tfrac{1}{4}%
+n\right) T\right) \right\Vert
\end{equation*}%
where $P_{K_{\gamma }}$ denotes the orthogonal projection of $L^{2}\left(
\Omega \right) $ onto the second subspace. Therefore, there is no interest
to control the initial data $P_{K_{\gamma }}y\left( \left( \frac{1}{4}%
+n\right) T\right) $ when $\lambda _{K_{\gamma }+1}\gg \gamma $. This
suggests the form of our feedback law.

In our analysis, a precise estimate is established and we will build the
output stabilization law via some minimal norm impulse control problems. It
requires to link the approximate impulse control and a quantitative unique
continuation estimate called observation at one time. Our program follows
the orientation described in \cite{Li} where stabilization, optimal control
and exact controllability for hyperbolic systems are closely linked. Here,
exact controllability for the wave equation is replaced by approximate
impulse control for the heat equations.

\bigskip

\bigskip

Several notes are given in order.

\begin{enumerate}
\item Impulse control belongs to a class of important control and has wide
applications. In many cases impulse control is an interesting alternative to
deal with systems that cannot be acted on by means of continuous control
inputs, for instance, relevant control for acting on a population of
bacteria should be impulsive, so that the density of the bactericide may
change instantaneously; indeed continuous control would be enhance drug
resistance of bacteria (see \cite{TWZ}). There are many studies on optimal
control and controllability for impulse controlled equations (see, for
instance, \cite{BL1}, \cite{R}, \cite{BL2}, \cite{LM}, \cite{Z}, \cite{LY}, 
\cite{BC}, \cite{DS}, \cite{MR}, \cite{OS}, \cite{Be} and references
therein). However, we have not found any published paper on stabilization
for impulse controlled equations. From this perspective, the problem studied
in the current paper is new.

In the system (\ref{1.1}), we do not need to follow the rule: Make
observation at each time and then add simultaneously control with a feedback
form. (When systems have continuous feedback control inputs, one has to
follow such a rule.) Instead of this, we only need to observe $1_{\omega
_{1}}^{\ast }y$ at time points $\left( n-\frac{1}{2}\right) T$ (for all $%
n\in \mathbb{N}$) and then add the controls at time points $nT$. In this
way, we can not only save observation and control time, but also allow the
control time $nT$ having a delay with respect to the observation time $%
\left( n-\frac{1}{2}\right) T$. Sometime, such time delay could be important
in practical application.

\item Stabilization is one of the most important subject in control theory.
In most studies of this subject, the aim of stabilization is to ask for a
feedback law so that the closed loop equation decays exponentially. The
current work aims to find, for each decay rate $\gamma $, a feedback law so
that the closed loop equation has an exponential decay with the rate $\gamma 
$. Such kind of stabilization is called the rapid stabilization. About this
subject, we would like to mention the works \cite{K}, \cite{U}, \cite{CCr}, 
\cite{V}, \cite{CCo}, \cite{CL2}, and \cite{CL1}.

\item When observation region is not the whole $\Omega $, the corresponding
stabilization is a kind of output stabilization. Such stabilization is very
useful in applications. Unfortunately, there is no systematic study on this
subject, even for the simplest case when the controlled system is
time-invariant linear ODE (see \cite{Br}). Most of publications on this
subject focus on how to construct an output feedback law for some special
equations (see, for instance, \cite{I}, \cite{Cu}, \cite{YY}, \cite{Co}, 
\cite{NS} and references therein). Our study also only provides an output
feedback law for a special equation.

\item In this paper, we present a new way to build up the feedback law. In
particular, the structure of our feedback law is not based on LQ theory or
Lyapunov functions (see, for instance, \cite{Ba} and \cite{C}).

\item In the current study, one of the keys to build up our feedback law is
the use of the unique continuation estimate at one time, established in \cite%
{PWZ} (see also \cite{PW} and \cite{PW1}). Some new observations are made on
it in this paper (see Theorem \ref{theorem2.1}, Remark \ref{remark2.2} and
Remark \ref{remark2.3}).

\item The following extensions of the current work should be interesting:
The first case is that $V$ depends on both $x$ and $t$ variables; The second
case is that the equation is a semi-linear heat equation; The third case is
that the equation is other types of PDEs.
\end{enumerate}

\bigskip

\bigskip

The rest of this paper is organized as follows. Section 2 provides several
inequalities which are equivalent to the unique continuation estimate at one
time. Such observation estimates are used in Section 3 in which we deal with
the impulse control problem. In Section 4, we link the impulse control
problem with an approximate inverse source problem. Finally, Section 5
presents the main result, as well as its proof.

\bigskip

\bigskip

\section{Observation at one time}

\bigskip

In this section, we present several equivalent inequalities. One of them is
the unique continuation estimate at one time built up in \cite{PWZ} (see
also \cite{PW} and \cite{PW1}).

\bigskip

\begin{thm}
\label{theorem2.1} Let $\omega $ be an open and nonempty subset of $\Omega $%
. Then the following propositions are equivalent and are true:

\begin{description}
\item[$\left( i\right) $] There are two constants $C_{1}>0$ and $\beta \in
\left( 0,1\right) $, which depend only on $\Omega $ and $\omega $, so that
for all $t>0$ and $\Phi \in L^{2}\left( \Omega \right) $,%
\begin{equation*}
\left\Vert e^{tA}\Phi \right\Vert \leq e^{C_{1}\left( 1+\frac{1}{t}%
+t\left\Vert V\right\Vert _{\infty }+\left\Vert V\right\Vert _{\infty
}^{2/3}\right) }\left\Vert \Phi \right\Vert ^{\beta }\left\Vert 1_{\omega
}^{\ast }e^{tA}\Phi \right\Vert _{\omega }^{1-\beta }\text{ .}
\end{equation*}

\item[$\left( ii\right) $] There is a positive constant $C_{2}$, depending
only on $\Omega $ and $\omega $, so that for each $\lambda \geq 0$ and each
sequence of real numbers $\left\{ a_{j}\right\} \subset \mathbb{R}$, 
\begin{equation*}
\sum_{\lambda _{j}<\lambda }\left\vert a_{j}\right\vert ^{2}\leq
e^{C_{2}\left( 1+\left\Vert V\right\Vert _{\infty }^{2/3}+\sqrt{\lambda }%
\right) }\int_{\omega }\left\vert \sum_{\lambda _{j}<\lambda }a_{j}\xi
_{j}\right\vert ^{2}dx\text{ .}
\end{equation*}

\item[$\left( iii\right) $] There is a positive constant $C_{3}$, depending
only on $\Omega $ and $\omega $, so that for all $\theta \in \left(
0,1\right) $, $t>0$ and $\Phi \in L^{2}\left( \Omega \right) $,%
\begin{equation*}
\left\Vert e^{tA}\Phi \right\Vert \leq e^{C_{3}\left( 1+\frac{1}{\theta t}%
+t\left\Vert V\right\Vert _{\infty }+\left\Vert V\right\Vert _{\infty
}^{2/3}\right) }\left\Vert \Phi \right\Vert ^{\theta }\left\Vert 1_{\omega
}^{\ast }e^{tA}\Phi \right\Vert _{\omega }^{1-\theta }\text{ .}
\end{equation*}

\item[$\left( iv\right) $] There is a positive constant $C_{3}$, depending
only on $\Omega $ and $\omega $, so that for all $\varepsilon ,\beta >0$, $%
t>0$ and $\Phi \in L^{2}\left( \Omega \right) $,%
\begin{equation*}
\left\Vert e^{tA}\Phi \right\Vert \leq \frac{1}{\varepsilon ^{\beta }}%
e^{C_{3}\left( 1+\beta \right) \left( 1+\frac{1+\beta }{\beta t}+t\left\Vert
V\right\Vert _{\infty }+\left\Vert V\right\Vert _{\infty }^{2/3}\right)
}\left\Vert 1_{\omega }^{\ast }e^{tA}\Phi \right\Vert _{\omega }+\varepsilon
\left\Vert \Phi \right\Vert \text{ .}
\end{equation*}

\item[$\left( v\right) $] There is a positive constant $c$, depending only
on $\Omega $ and $\omega $, so that for all $\varepsilon >0$, $t>0$ and $%
\Phi \in L^{2}\left( \Omega \right) $,%
\begin{equation*}
\left\Vert e^{tA}\Phi \right\Vert \leq e^{c\left( 1+\frac{1}{t}+t\left\Vert
V\right\Vert _{\infty }+\left\Vert V\right\Vert _{\infty }^{2/3}\right) }%
\normalfont{\text{exp}}\left( \sqrt{\frac{c}{t}\normalfont{\text{ln}}^{+}%
\frac{1}{\varepsilon }}\right) \left\Vert 1_{\omega }^{\ast }e^{tA}\Phi
\right\Vert _{\omega }+\varepsilon \left\Vert \Phi \right\Vert \text{ .}
\end{equation*}
\end{description}

Here $\normalfont{\text{ln}}^{+}\frac{1}{\varepsilon }:=\normalfont{%
\text{max}}\{\normalfont{\text{ln}}\frac{1}{\varepsilon },0\}$. Moreover,
constants $C_{3}$ in $\left( iii\right) $ and $\left( iv\right) $ can be
chosen as the same number, and $c=4C_{3}$ in $\left( v\right) $.
\end{thm}

\bigskip

\begin{rem}
\label{remark2.2} We would like to give several notes on Theorem~\ref%
{theorem2.1}.

\begin{enumerate}
\item The inequality in $(i)$ of Theorem \ref{theorem2.1} implies a
quantitative version of the unique continuation property of heat equations
built up in \cite{L} (see also \cite{EFV}). Indeed, if $\left\Vert 1_{\omega
}^{\ast }e^{tA}\Phi \right\Vert _{\omega }=0$ for some $\Phi \in
L^{2}(\Omega )$, then $\left\Vert e^{tA}\Phi \right\Vert =0$, which,
together with the backward uniqueness for heat equations, implies that $\Phi
=0$. Moreover, such inequality is equivalent to a kind of impulse control,
with a bound on its cost, which will be explained in the next section.
Further, such controllability is the base for our study on the stabilization.

\item In the following studies of this paper, we will essentially use the
inequality in $(iii)$ of Theorem \ref{theorem2.1}. However, other
inequalities seem to be interesting independently. For instance, when $V=0$,
the inequality in $(ii)$ of Theorem \ref{theorem2.1} is exactly the
Lebeau-Robbiano spectral inequality (see \cite{LR}, \cite{JL} or \cite{LZ}).
Here, we get, in the case that $V\neq 0$, the same inequality, and find how
the constant (on the right hand side of the inequality) depends on $%
\left\Vert V\right\Vert _{\infty }$.

\item It deserves to mention that all constant terms on the right hand side
of inequalities in Theorem \ref{theorem2.1} have explicit expressions in
terms of the norm of $V$ and time, but not $\Omega $ and $\omega $.

\item It was realized that when $V=0$, the inequality in $(i)$ of Theorem %
\ref{theorem2.1} can imply the inequality $(ii)$ of Theorem \ref{theorem2.1}
(see Remark 1 in \cite{AEWZ}).

\item The key ingredient why the inequality in $(i)$ and the one in $(ii)$
of Theorem \ref{theorem2.1} are equivalent is that the evolution $e^{tA}$ is
time-invariant. We would like to mention that the inequality in $(i)$ of
Theorem \ref{theorem2.1} was proved in \cite{PWZ} (see also \cite{PW} and 
\cite{PW1}) for the case when $V=V(x,t)$ (i.e., for the evolutions with
time-varying coefficients). However, for the time-varying evolutions, we do
not know whether $(iii)$, $(iv)$ and $(v)$ in Theorem \ref{theorem2.1} are
still valid and what should be the right alternative of the spectral
inequality in $(ii)$ of Theorem \ref{theorem2.1}. These should be
interesting open problems.
\end{enumerate}
\end{rem}

\bigskip

\bigskip

\begin{rem}
\label{remark2.3} When $V=0$, the function $\normalfont{\text{exp}}\left( 
\sqrt{\normalfont{\text{ln}}^{+}\frac{1}{\varepsilon }}\right) $ is optimal
in the inequality in $(v)$ of Theorem \ref{theorem2.1} in the following
sense: For any function $\varepsilon \mapsto f\left( \varepsilon \right) $, $%
\varepsilon >0$, with $\underset{\varepsilon \rightarrow 0^{+}}{\overline{%
\normalfont{\text{lim}}}}\frac{f\left( \varepsilon \right) }{\sqrt{%
\normalfont{\text{ln}}^{+}\frac{1}{\varepsilon }}}=0$, the following
inequality is not true:%
\begin{equation*}
\left\Vert e^{t\Delta }\Phi \right\Vert \leq f\left( \varepsilon \right)
e^{f\left( \varepsilon \right) }\left\Vert 1_{\omega }^{\ast }e^{t\Delta
}\Phi \right\Vert _{\omega }+\varepsilon \left\Vert \Phi \right\Vert \quad 
\text{for all }\Phi \in L^{2}(\Omega )\text{ and }\varepsilon >0\text{ .}
\end{equation*}
\end{rem}

\bigskip

We now explain why the optimality in Remark~\ref{remark2.3} is true.
According to Proposition 5.5 in \cite{LL}, there exist $C_{0}^{\prime }>0$
and $n_{0}>0$ so that for each $n\geq n_{0}$, there is $\left\{
a_{n,j}\right\} _{j=1}^{\infty }\subset \ell ^{2}\left\backslash \left\{
0\right\} \right. $ satisfying that%
\begin{equation*}
\left\Vert \sum_{\lambda _{j}\leq n}a_{n,j}\xi _{j}\right\Vert \geq
C_{0}^{\prime }e^{C_{0}^{\prime }\sqrt{n}}\left\Vert 1_{\omega }^{\ast
}\sum_{\lambda _{j}\leq n}a_{n,j}\xi _{j}\right\Vert _{\omega }\text{ .}
\end{equation*}%
Define, for each $n\in \mathbb{N}$,%
\begin{equation*}
\Phi _{n}:=\sum_{\lambda _{j}\leq n}e^{\lambda _{j}t}a_{n,j}\xi _{j}\text{
and }\delta _{n}:=\frac{1}{2}e^{-nt}\text{ .}
\end{equation*}%
Therefore, we find that for each $n\geq n_{0}$,%
\begin{equation*}
\begin{array}{ll}
\left\Vert e^{t\Delta }\Phi _{n}\right\Vert & =\left\Vert \displaystyle%
\sum_{\lambda _{j}\leq n}a_{n,j}\xi _{j}\right\Vert =\left( 1-\delta
_{n}e^{nt}\right) \left\Vert \displaystyle\sum_{\lambda _{j}\leq
n}a_{n,j}\xi _{j}\right\Vert +\delta _{n}e^{nt}\left\Vert \displaystyle%
\sum_{\lambda _{j}\leq n}a_{n,j}\xi _{j}\right\Vert \\ 
& \geq \displaystyle\frac{1}{2}\left\Vert \displaystyle\sum_{\lambda
_{j}\leq n}a_{n,j}\xi _{j}\right\Vert +\delta _{n}\left\Vert \displaystyle%
\sum_{\lambda _{j}\leq n}e^{\lambda _{j}t}a_{n,j}\xi _{j}\right\Vert \\ 
& \geq \displaystyle\frac{1}{2}C_{0}^{\prime }e^{C_{0}^{\prime }\sqrt{n}%
}\left\Vert 1_{\omega }^{\ast }\displaystyle\sum_{\lambda _{j}\leq
n}a_{n,j}\xi _{j}\right\Vert _{\omega }+\delta _{n}\left\Vert \Phi
_{n}\right\Vert \text{ .}%
\end{array}%
\end{equation*}%
Meanwhile, it follows from the definition of $\delta _{n}$ that for each $%
n\in \mathbb{N}$,%
\begin{equation*}
\sqrt{\text{ln}\frac{1}{\delta _{n}}}=\sqrt{\text{ln}2+nt}\leq 1+\sqrt{n}%
\sqrt{t}\text{ .}
\end{equation*}%
Gathering all the previous estimates, we see that%
\begin{equation*}
\left\Vert e^{t\Delta }\Phi _{n}\right\Vert \geq \frac{1}{2}C_{0}^{\prime
}e^{-\frac{C_{0}^{\prime }}{\sqrt{t}}}e^{\frac{C_{0}^{\prime }}{\sqrt{t}}%
\sqrt{\text{ln}\frac{1}{\delta _{n}}}}\left\Vert 1_{\omega }^{\ast
}\sum_{\lambda _{j}\leq n}a_{n,j}\xi _{j}\right\Vert _{\omega }+\delta
_{n}\left\Vert \Phi _{n}\right\Vert \text{ ,}
\end{equation*}%
which leads to the following property: There exists $C_{0}>0$, $\left\{ \Phi
_{n}\right\} \subset L^{2}\left( \Omega \right) \left\backslash \left\{
0\right\} \right. $ and $\left\{ \delta _{n}\right\} \subset \left(
0,1\right) $, with $\underset{n\rightarrow \infty }{\text{lim}}\delta _{n}=0$%
, so that%
\begin{equation*}
\left\Vert e^{t\Delta }\Phi _{n}\right\Vert \geq C_{0}e^{C_{0}\sqrt{\text{ln}%
\frac{1}{\delta _{n}}}}\left\Vert 1_{\omega }^{\ast }e^{t\Delta }\Phi
_{n}\right\Vert _{\omega }+\delta _{n}\left\Vert \Phi _{n}\right\Vert \quad 
\text{for all }n\text{ .}
\end{equation*}

\bigskip

The rest of this section is devoted to the proof of Theorem \ref{theorem2.1}.

\bigskip

Proof. We organize the proof of Theorem \ref{theorem2.1} by several steps.

\noindent\textit{Step 1: On the proposition }$(i)$.

The conclusion $(i)$ has been proved in \cite{PWZ} (see also \cite{PW} and 
\cite{PW1}).

\noindent\textit{Step 2: To prove that }$(i)\Rightarrow(ii)$.

Let $C_{1}>0$ and $\beta \in (0,1)$ be given by $(i)$. Arbitrarily fix $%
\lambda \geq 0$ and $\{a_{j}\}\subset \mathbb{R}$. By applying the
inequality in $(i)$, with $\Phi =\sum_{\lambda _{j}<\lambda }a_{j}e^{\lambda
_{j}t}\xi _{j}$, we get that 
\begin{equation*}
\sum_{\lambda _{j}<\lambda }\left\vert a_{j}\right\vert ^{2}\leq
e^{2C_{1}\left( 1+\frac{1}{t}+t\left\Vert V\right\Vert _{\infty }+\left\Vert
V\right\Vert _{\infty }^{2/3}\right) }\left( \sum_{\lambda _{j}<\lambda
}\left\vert a_{j}e^{\lambda _{j}t}\right\vert ^{2}\right) ^{\beta }\left(
\int_{\omega }\left\vert \sum_{\lambda _{j}<\lambda }a_{j}\xi
_{j}\right\vert ^{2}dx\right) ^{1-\beta }\text{ ,}
\end{equation*}%
which implies that 
\begin{equation}
\sum_{\lambda _{j}<\lambda }|a_{j}|^{2}\leq e^{\frac{2}{1-\beta }C_{1}\left(
1+\frac{1}{t}+t\Vert V\Vert _{\infty }+\Vert V\Vert _{\infty }^{2/3}\right)
}e^{\frac{2\beta }{1-\beta }\lambda t}\int_{\omega }\left\vert \sum_{\lambda
_{j}<\lambda }a_{j}\xi _{j}\right\vert ^{2}dx\quad \text{\ for each }t>0%
\text{ .}  \tag{2.1}  \label{2.1}
\end{equation}%
Meanwhile, since $\Vert V\Vert _{\infty }^{1/2}\leq 1+\Vert V\Vert _{\infty
}^{2/3}$ and $\beta \in (0,1)$, we see that%
\begin{equation*}
\begin{array}{ll}
& \quad \underset{t>0}{\text{inf}}\left[ \displaystyle C_{1}\left( 1+\frac{1%
}{t}+t\Vert V\Vert _{\infty }+\Vert V\Vert _{\infty }^{2/3}\right) +\beta
\lambda t\right] \\ 
& =\underset{t>0}{\text{inf}}\left[ \displaystyle C_{1}(1+\Vert V\Vert
_{\infty }^{2/3})+\frac{C_{1}}{t}+\left( C_{1}\Vert V\Vert _{\infty }+\beta
\lambda \right) t\right] \\ 
& =C_{1}\left( 1+\Vert V\Vert _{\infty }^{2/3}\right) +2\sqrt{C_{1}\left(
C_{1}\Vert V\Vert _{\infty }+\beta \lambda \right) } \\ 
& \leq C_{1}\left( 1+\Vert V\Vert _{\infty }^{2/3}+2\Vert V\Vert _{\infty
}^{1/2}\right) +2\sqrt{C_{1}}\sqrt{\beta \lambda }\leq \text{max}\left\{
3C_{1},2\sqrt{C_{1}}\right\} \left( 1+\Vert V\Vert _{\infty }^{2/3}+\sqrt{%
\lambda }\right) \text{ .}%
\end{array}%
\end{equation*}%
This, along with (\ref{2.1}), leads to the conclusion $(ii)$, with $C_{2}=$%
max$\left\{ \frac{6C_{1}}{1-\beta },\frac{4\sqrt{C_{1}}}{1-\beta }\right\} $.

\noindent\textit{Step 3: To show that }$(ii)\Rightarrow(iii)$.

Arbitrarily fix $\lambda \geq 0$, $t>0$ and $\Phi =\sum_{j\geq 1}a_{j}\xi
_{j}$ with $\{a_{j}\}\subset \ell ^{2}$. Write 
\begin{equation*}
e^{tA}\Phi =\sum_{\lambda _{j}<\lambda }a_{j}e^{-\lambda _{j}t}\xi
_{j}+\sum_{\lambda _{j}\geq \lambda }a_{j}e^{-\lambda _{j}t}\xi _{j}\text{ .}
\end{equation*}%
Then by $(ii)$, we find that%
\begin{equation*}
\begin{array}{ll}
\left\Vert e^{tA}\Phi \right\Vert & \leq \left\Vert \displaystyle%
\sum_{\lambda _{j}<\lambda }a_{j}e^{-\lambda _{j}t}\xi _{j}\right\Vert
+\left\Vert \displaystyle\sum_{\lambda _{j}\geq \lambda }a_{j}e^{-\lambda
_{j}t}\xi _{j}\right\Vert \\ 
& \leq \left( \displaystyle\sum_{\lambda _{j}<\lambda }\left\vert
a_{j}e^{-\lambda _{j}t}\right\vert ^{2}\right) ^{1/2}+e^{-\lambda
t}\left\Vert \Phi \right\Vert \\ 
& \leq \left( e^{C_{2}\left( 1+\Vert V\Vert _{\infty }^{2/3}+\sqrt{\lambda }%
\right) }\displaystyle\int_{\omega }\left\vert \displaystyle\sum_{\lambda
_{j}<\lambda }a_{j}e^{-\lambda _{j}t}\xi _{j}\right\vert ^{2}dx\right)
^{1/2}+e^{-\lambda t}\left\Vert \Phi \right\Vert \text{ .}%
\end{array}%
\end{equation*}%
This, along with the triangle inequality for the norm $\left\Vert \cdot
\right\Vert _{\omega }$, yields that%
\begin{equation*}
\begin{array}{ll}
\left\Vert e^{tA}\Phi \right\Vert & \leq \left( e^{C_{2}\left( 1+\Vert
V\Vert _{\infty }^{2/3}+\sqrt{\lambda }\right) }\displaystyle\int_{\omega
}\left\vert \displaystyle\sum_{j\geq 1}a_{j}e^{-\lambda _{j}t}\xi
_{j}\right\vert ^{2}dx\right) ^{1/2} \\ 
& \quad +\left( e^{C_{2}\left( 1+\Vert V\Vert _{\infty }^{2/3}+\sqrt{\lambda 
}\right) }\displaystyle\int_{\omega }\left\vert \displaystyle\sum_{\lambda
_{j}\geq \lambda }a_{j}e^{-\lambda _{j}t}\xi _{j}\right\vert ^{2}dx\right)
^{1/2}+e^{-\lambda t}\left\Vert \Phi \right\Vert \text{ .}%
\end{array}%
\end{equation*}%
Hence, it follows that%
\begin{equation*}
\begin{array}{ll}
\left\Vert e^{tA}\Phi \right\Vert & \leq e^{\frac{C_{2}}{2}\left( 1+\Vert
V\Vert _{\infty }^{2/3}+\sqrt{\lambda }\right) }\left\Vert 1_{\omega }^{\ast
}e^{tA}\Phi \right\Vert _{\omega }+e^{\frac{C_{2}}{2}\left( 1+\Vert V\Vert
_{\infty }^{2/3}+\sqrt{\lambda }\right) }e^{-\lambda t}\left\Vert \Phi
\right\Vert +e^{-\lambda t}\left\Vert \Phi \right\Vert \\ 
& \leq 2e^{\frac{C_{2}}{2}\left( 1+\Vert V\Vert _{\infty }^{2/3}\right) }e^{%
\frac{C_{2}}{2}\sqrt{\lambda }}\left( \left\Vert 1_{\omega }^{\ast
}e^{tA}\Phi \right\Vert _{\omega }+e^{-\lambda t}\left\Vert \Phi \right\Vert
\right) \text{ .}%
\end{array}%
\end{equation*}%
Combining the above estimate with the following inequality: 
\begin{equation*}
\frac{C_{2}}{2}\sqrt{\lambda }\leq \frac{\rho }{2}\lambda t+\frac{1}{2t\rho }%
\left( \frac{C_{2}}{2}\right) ^{2}\quad \text{for any }\rho >0\text{ ,}
\end{equation*}%
we have that for all $\rho \in (0,2)$ and $\lambda \geq 0$, 
\begin{equation*}
\left\Vert e^{tA}\Phi \right\Vert \leq 2e^{\frac{C_{2}}{2}\left( 1+\Vert
V\Vert _{\infty }^{2/3}\right) }e^{\frac{1}{2t\rho }\left( \frac{C_{2}}{2}%
\right) ^{2}}\left( e^{\frac{\rho }{2}\lambda t}\left\Vert 1_{\omega }^{\ast
}e^{tA}\Phi \right\Vert _{\omega }+e^{-\frac{2-\rho }{2}\lambda t}\left\Vert
\Phi \right\Vert \right) \text{ .}
\end{equation*}%
Since $\lambda $ was arbitrarily taken from $[0,\infty )$, we choose 
\begin{equation*}
\lambda =\frac{1}{t}\text{ln}\left( \frac{e^{t\Vert V\Vert _{\infty
}}\left\Vert \Phi \right\Vert }{\left\Vert 1_{\omega }^{\ast }e^{tA}\Phi
\right\Vert _{\omega }}\right)
\end{equation*}%
to get 
\begin{equation*}
\left\Vert e^{tA}\Phi \right\Vert \leq 2e^{\frac{C_{2}}{2}\left( 1+\Vert
V\Vert _{\infty }^{2/3}\right) }e^{\frac{1}{2t\rho }\left( \frac{C_{2}}{2}%
\right) ^{2}}\left( 2e^{t\Vert V\Vert _{\infty }}\left\Vert 1_{\omega
}^{\ast }e^{tA}\Phi \right\Vert _{\omega }^{1-\frac{\rho }{2}}\left\Vert
\Phi \right\Vert ^{\frac{\rho }{2}}\right)
\end{equation*}%
which is the inequality in $(iii)$ with $\theta =\rho /2$.

\noindent \textit{Step 4: To show that }$(iii)\Rightarrow (iv)$.

We write the inequality in $(iii)$ in the following way: 
\begin{equation*}
\left\Vert e^{tA}\Phi \right\Vert \leq \left\Vert \Phi \right\Vert ^{\theta
}\left( e^{\frac{C_{3}}{1-\theta }\left( 1+\frac{1}{\theta t}+t\left\Vert
V\right\Vert _{\infty }+\left\Vert V\right\Vert _{\infty }^{2/3}\right)
}\left\Vert 1_{\omega }^{\ast }e^{tA}\Phi \right\Vert _{\omega }\right)
^{1-\theta }\text{ .}
\end{equation*}%
Notice that for any real numbers $E,B,D>0$ and $\theta \in (0,1)$ 
\begin{equation*}
\begin{array}{ll}
& E\leq B^{\theta }D^{1-\theta }\text{ } \\ 
\Leftrightarrow & E\leq \varepsilon B+\left( 1-\theta \right) \theta ^{\frac{%
\theta }{1-\theta }}\displaystyle\frac{1}{\varepsilon ^{\frac{\theta }{%
1-\theta }}}D\quad \forall \varepsilon >0\text{ .}%
\end{array}%
\end{equation*}%
Then by taking $\beta =\frac{\theta }{1-\theta }$ in the last inequality, we
are led to the inequality in $(iv)$ with the same constant $C_{3}$ as that
in $(iii)$. (The above equivalence can be easily verified by using the Young
inequality and by choosing $\varepsilon =\theta \left( \frac{D}{B}\right)
^{1-\theta }$.)

\noindent \textit{Step 5: To show that }$(iv)\Rightarrow (v)$.

We write the inequality in $(iv)$ in the following way:%
\begin{equation*}
\left\Vert e^{tA}\Phi \right\Vert \leq \text{exp}\left( \Lambda +\Upsilon
+\beta \left( \text{ln}\frac{1}{\varepsilon }+\Lambda \right) +\frac{1}{%
\beta }\Upsilon \right) \left\Vert 1_{\omega }^{\ast }e^{tA}\Phi \right\Vert
_{\omega }+\varepsilon \left\Vert \Phi \right\Vert \text{ ,}
\end{equation*}%
with $\Lambda =C_{3}\left( 1+\frac{1}{t}+t\left\Vert V\right\Vert _{\infty
}+\left\Vert V\right\Vert _{\infty }^{2/3}\right) $ and $\Upsilon =\frac{%
C_{3}}{t}$. Next, we optimize the above inequality with respect to $\beta >0$
by choosing $\beta =\sqrt{\frac{\Upsilon }{\text{ln}^{+}\frac{1}{\varepsilon 
}+\Lambda }}$ to get%
\begin{equation*}
\begin{array}{ll}
\left\Vert e^{tA}\Phi \right\Vert & \leq \text{exp}\left( \Lambda +\Upsilon
+2\sqrt{\Upsilon \left( \text{ln}^{+}\frac{1}{\varepsilon }+\Lambda \right) }%
\right) \left\Vert 1_{\omega }^{\ast }e^{tA}\Phi \right\Vert _{\omega
}+\varepsilon \left\Vert \Phi \right\Vert \\ 
& \leq \text{exp}\left( \Lambda +\Upsilon +2\sqrt{\Upsilon \Lambda }+2\sqrt{%
\Upsilon \text{ln}^{+}\frac{1}{\varepsilon }}\right) \left\Vert 1_{\omega
}^{\ast }e^{tA}\Phi \right\Vert _{\omega }+\varepsilon \left\Vert \Phi
\right\Vert \\ 
& \leq \text{exp}\left( 4\Lambda +2\sqrt{\Upsilon \text{ln}^{+}\frac{1}{%
\varepsilon }}\right) \left\Vert 1_{\omega }^{\ast }e^{tA}\Phi \right\Vert
_{\omega }+\varepsilon \left\Vert \Phi \right\Vert \text{ .}%
\end{array}%
\end{equation*}%
This implies the inequality in $(v)$, with $c=4C_{3}$.

\noindent \textit{Step 6: to show that }$(v)\Rightarrow (i)$.

Since 
\begin{equation*}
\sqrt{\frac{c}{t}\text{ln}^{+}\frac{1}{\varepsilon }}\leq \frac{c}{\alpha t}%
+\alpha \text{ln}\left( e+\frac{1}{\varepsilon }\right) \quad \forall \alpha
>0\text{ ,}
\end{equation*}%
the inequality in $(v)$ becomes%
\begin{equation*}
\left\Vert e^{tA}\Phi \right\Vert \leq e^{c\left( 1+\frac{1}{t}+t\left\Vert
V\right\Vert _{\infty }+\left\Vert V\right\Vert _{\infty }^{2/3}\right) }e^{%
\frac{c}{\alpha t}}\left( e+\frac{1}{\varepsilon }\right) ^{\alpha
}\left\Vert 1_{\omega }^{\ast }e^{tA}\Phi \right\Vert _{\omega }+\varepsilon
\left\Vert \Phi \right\Vert \text{ .}
\end{equation*}%
Next, we choose 
\begin{equation*}
\varepsilon =\frac{1}{2}\frac{\left\Vert e^{tA}\Phi \right\Vert }{\left\Vert
\Phi \right\Vert }
\end{equation*}%
and we use the fact that $\left\Vert e^{tA}\Phi \right\Vert \leq e^{t\Vert
V\Vert _{\infty }}\left\Vert \Phi \right\Vert $ to deduce the inequality in $%
(i)$.

This ends the proof.

\bigskip

\bigskip

\section{Impulse control}

\bigskip

In this section, we first state our key result on impulse approximate
controllability.

\bigskip

\begin{thm}
\label{theorem3.1} Let $0\leq T_{1}<T_{2}<T_{3}$. Let $\varepsilon >0$ and $%
z\in L^{2}\left( \Omega \right) $. Then the following conclusions are true:

\begin{description}
\item[$\left( i\right) $] There exists $f\in L^{2}\left( \omega \right) $\
such that the unique solution $y$ to the equation: 
\begin{equation*}
\left\{ 
\begin{array}{ll}
y^{\prime }\left( t\right) -Ay\left( t\right) =0\text{ ,} & t\in \left(
T_{1},T_{3}\right) \backslash \left\{ T_{2}\right\} \ \text{,} \\ 
y\left( T_{1}\right) =z\text{ ,} &  \\ 
y\left( T_{2}\right) =y\left( T_{2-}\right) +1_{\omega }f\text{ ,} & 
\end{array}%
\right.
\end{equation*}%
satisfies 
\begin{equation*}
\left\Vert y\left( T_{3}\right) \right\Vert \leq \varepsilon \left\Vert
z\right\Vert \text{ .}
\end{equation*}%
Moreover, it holds that 
\begin{equation*}
\left\Vert f\right\Vert _{\omega }\leq \mathcal{C}_{\varepsilon }\left(
T_{3}-T_{2},T_{2}-T_{1}\right) \left\Vert z\right\Vert
\end{equation*}%
where $\mathcal{C}_{\varepsilon }:\mathbb{R}^{+}\times \mathbb{R}%
^{+}\rightarrow \mathbb{R}^{+}$ is given by 
\begin{equation*}
\mathcal{C}_{\varepsilon }\left( t,s\right) =e^{4s\left\Vert V\right\Vert
_{\infty }}e^{c\left( 1+\frac{1}{t}+t\left\Vert V\right\Vert _{\infty
}+\left\Vert V\right\Vert _{\infty }^{2/3}\right) }\normalfont{\text{exp}}%
\left( \sqrt{\frac{c}{t}\normalfont{\text{ln}}^{+}\frac{1}{\varepsilon }}%
\right) \text{ ,}
\end{equation*}%
and $c>0$ is the same constant than in $(v)$ of Theorem 2.1 .

\item[$\left( ii\right) $] There is a unique $f$ solving the following
problem: 
\begin{equation*}
\normalfont{\text{inf}}\left\{ \left\Vert h\right\Vert _{\omega };h\in
L^{2}\left( \omega \right) \text{ such that }\left\Vert y\left( T_{3}\right)
\right\Vert \leq \varepsilon \left\Vert z\right\Vert \text{ with }f\text{
replaced by }h\text{ in }\left( i\right) \right\} \text{ .}
\end{equation*}
\end{description}
\end{thm}

\bigskip

\bigskip

The proof of Theorem \ref{theorem3.1} will be given in Subsections 3.1-3.2.
We now apply Theorem \ref{theorem3.1} to the eigenfunctions $\left\{ \xi
_{j}\right\} _{j=1}^{\infty }$. More precisely, by choosing $z=\xi _{j}$ in
Theorem \ref{theorem3.1}, we have the following corollary:

\begin{cor}
\label{corollary3.2} For any $0\leq T_{1}<T_{2}<T_{3}$ and any $\varepsilon
>0$, $j\in \mathbb{N}$, there is a pair $\left( y_{j},f_{j}\right) $ such
that 
\begin{equation*}
\left\{ 
\begin{array}{ll}
y_{j}^{\prime }\left( t\right) -Ay_{j}\left( t\right) =0\text{ ,} & t\in
\left( T_{1},T_{3}\right) \backslash \left\{ T_{2}\right\} \ \text{,} \\ 
y_{j}\left( T_{1}\right) =\xi _{j}\text{ ,} &  \\ 
y_{j}\left( T_{2}\right) =y_{j}\left( T_{2-}\right) +1_{\omega }f_{j}\text{ ,%
} & 
\end{array}%
\right.
\end{equation*}%
and 
\begin{equation*}
\left\{ 
\begin{array}{ll}
\left\Vert y_{j}\left( T_{3}\right) \right\Vert \leq \varepsilon \text{ ,} & 
\\ 
\left\Vert f_{j}\right\Vert _{\omega }\leq \mathcal{C}_{\varepsilon }\left(
T_{3}-T_{2},T_{2}-T_{1}\right) \text{ ,} & 
\end{array}%
\right.
\end{equation*}%
where $\mathcal{C}_{\varepsilon }$ is given in Theorem \ref{theorem3.1}.
Further, the control function $f_{j}$\ can be taken as the unique solution
of the problem 
\begin{equation*}
\normalfont{\text{inf}}\left\{ 
\begin{array}{c}
\left\Vert h\right\Vert _{\omega };h\in L^{2}\left( \omega \right) \text{ }%
and\text{ }the\text{ }property\text{ }\left( i\right) \text{ }of\text{ }%
Theorem\text{ }\ref{theorem3.1} \\ 
holds\text{ }with\text{ }\left( z,f\right) \text{ }replaced\text{ }by\text{ }%
\left( \xi _{j},h\right)%
\end{array}%
\right\} \text{ .}
\end{equation*}
\end{cor}

\bigskip

\bigskip

Further we will apply Theorem \ref{theorem3.1} to a finite combination of
eigenfunctions. More precisely, we have the following consequence:

\bigskip

\begin{thm}
\label{theorem3.3} Let $\omega _{2}$ be a non-empty open subset of $\Omega $
and $K\in \mathbb{N}$. Let $0\leq T_{1}<T_{2}<T_{3}$ and $\varepsilon >0$.
Then for any $b=\left( b_{j}\right) _{j=1,\cdot \cdot ,K}$, there is a pair $%
\left( \widetilde{y},\widetilde{f}\right) $ such that 
\begin{equation*}
\left\{ 
\begin{array}{ll}
\widetilde{y}^{\prime }\left( t\right) -A\widetilde{y}\left( t\right) =0%
\text{ ,} & t\in \left( T_{1},T_{3}\right) \backslash \left\{ T_{2}\right\}
\ \text{,} \\ 
\widetilde{y}\left( T_{1}\right) =\displaystyle\sum_{j=1,\cdot \cdot
,K}b_{j}\xi _{j}\text{ ,} &  \\ 
\widetilde{y}\left( T_{2}\right) =\widetilde{y}\left( T_{2-}\right)
+1_{\omega _{2}}\widetilde{f}\text{ ,} & 
\end{array}%
\right.
\end{equation*}%
and 
\begin{equation*}
\left\{ 
\begin{array}{ll}
\left\Vert \widetilde{y}\left( T_{3}\right) \right\Vert \leq \varepsilon 
\sqrt{K}\left\Vert b\right\Vert _{\ell ^{2}}\text{ ,} &  \\ 
\widetilde{f}=\displaystyle\sum_{j=1,\cdot \cdot ,K}b_{j}f_{j}\text{ ,} & 
\end{array}%
\right.
\end{equation*}%
where $f_{j}$ is given by Corollary \ref{corollary3.2} with $\omega =\omega
_{2}$ and satisfies%
\begin{equation*}
\left\Vert f_{j}\right\Vert _{\omega _{2}}\leq \mathcal{C}_{\varepsilon
}\left( T_{3}-T_{2},T_{2}-T_{1}\right) \text{ .}
\end{equation*}
\end{thm}

\bigskip

In the study of our stabilization, we will use Corollary \ref{corollary3.2}
and Theorem \ref{theorem3.3}. The rest of this section is devoted to the
proof of Theorem \ref{theorem3.1}, and the studies on some minimal norm
control problem.

\bigskip

\subsection{Existence of impulse control functions and its cost}

\bigskip

The aim of this subsection is to prove the conclusion $(i)$ of Theorem \ref%
{theorem3.1}. It deserves mentioning what follows: The existence of controls
with $(i)$ of Theorem \ref{theorem3.1} is indeed the existence of impulse
controls (with a cost) driving the solution of the equation in $\left(
i\right) $ of Theorem \ref{theorem3.1} from the initial state $z$ to the
closed ball in $L^{2}(\Omega )$, centered at the origin and of radius $%
\varepsilon \left\Vert z\right\Vert $, at the ending time $T_{3}$.

\bigskip

To prove the conclusion $(i)$ of Theorem \ref{theorem3.1}, we let $%
\varepsilon >0$ and $z\in L^{2}(\Omega )$. Denote $\hbar =\varepsilon ^{2}$
and $k=\left( \mathcal{C}_{\varepsilon }\left(
T_{3}-T_{2},T_{2}-T_{1}\right) \right) ^{2}$. Consider the strictly convex $%
C^{1}$ functional $\digamma $ defined on $L^{2}\left( \Omega \right) $ given
by%
\begin{equation*}
\digamma \left( \Phi \right) :=\frac{k}{2}\left\Vert 1_{\omega }^{\ast
}e^{\left( T_{3}-T_{2}\right) A}\Phi \right\Vert _{\omega }^{2}+\frac{\hbar 
}{2}\left\Vert \Phi \right\Vert ^{2}+\left\langle z,e^{\left(
T_{3}-T_{1}\right) A}\Phi \right\rangle \text{ .}
\end{equation*}%
Notice that $\digamma $ is coercive and therefore $\digamma $ has a unique
minimizer $w\in L^{2}\left( \Omega \right) $, i.e. $\digamma (w)=\underset{%
\Phi \in L^{2}\left( \Omega \right) }{\text{min}}\digamma (\Phi )$. Since $%
\digamma ^{\prime }(w)\Phi =0$ for any $\Phi \in L^{2}\left( \Omega \right) $%
, we have%
\begin{equation}
k\left\langle 1_{\omega }^{\ast }e^{\left( T_{3}-T_{2}\right) A}w,1_{\omega
}^{\ast }e^{\left( T_{3}-T_{2}\right) A}\Phi \right\rangle _{\omega }+\hbar
\left\langle w,\Phi \right\rangle -\left\langle z,e^{\left(
T_{3}-T_{1}\right) A}\Phi \right\rangle =0\quad \forall \Phi \in L^{2}\left(
\Omega \right) \text{ .}  \tag{3.1.1}  \label{3.1.1}
\end{equation}%
But by multiplying by $e^{\left( T_{3}-t\right) A}\Phi $ the system solved
by $y$, one gets%
\begin{equation*}
\left\langle y\left( T_{3}\right) ,\Phi \right\rangle =\left\langle
z,e^{\left( T_{3}-T_{1}\right) A}\Phi \right\rangle +\left\langle 1_{\omega
}f,e^{\left( T_{3}-T_{2}\right) A}\Phi \right\rangle \quad \forall \Phi \in
L^{2}\left( \Omega \right) \text{ .}
\end{equation*}%
By choosing $f=-k1_{\omega }^{\ast }e^{\left( T_{3}-T_{2}\right) A}w$, the
above two equalities yield that 
\begin{equation*}
y\left( T_{3}\right) =\hbar w\text{ .}
\end{equation*}%
Further, one can deduce that 
\begin{equation*}
\frac{1}{k}\left\Vert f\right\Vert _{\omega }^{2}+\frac{1}{\hbar }\left\Vert
y\left( T_{3}\right) \right\Vert ^{2}=k\left\Vert 1_{\omega }^{\ast
}e^{\left( T_{3}-T_{2}\right) A}w\right\Vert _{\omega }^{2}+\hbar \left\Vert
w\right\Vert ^{2}\text{ .}
\end{equation*}%
Next, notice that by taking $\Phi =w$ in (\ref{3.1.1}) and by Cauchy-Schwarz
inequality, we have%
\begin{equation*}
k\left\Vert 1_{\omega }^{\ast }e^{\left( T_{3}-T_{2}\right) A}w\right\Vert
_{\omega }^{2}+\hbar \left\Vert w\right\Vert ^{2}=\left\langle z,e^{\left(
T_{3}-T_{1}\right) A}w\right\rangle \leq \frac{1}{2}\left\Vert z\right\Vert
^{2}+\frac{1}{2}\left\Vert e^{\left( T_{3}-T_{1}\right) A}w\right\Vert ^{2}%
\text{ .}
\end{equation*}%
Now, we claim that 
\begin{equation*}
\left\Vert e^{\left( T_{3}-T_{1}\right) A}w\right\Vert ^{2}\leq k\left\Vert
1_{\omega }^{\ast }e^{\left( T_{3}-T_{2}\right) A}w\right\Vert _{\omega
}^{2}+\hbar \left\Vert w\right\Vert ^{2}\text{ .}
\end{equation*}%
When it is proved, we can gather the previous three estimates to yield%
\begin{equation*}
\frac{1}{k}\left\Vert f\right\Vert _{\omega }^{2}+\frac{1}{\hbar }\left\Vert
y\left( T_{3}\right) \right\Vert ^{2}\leq \left\Vert z\right\Vert ^{2}\text{
.}
\end{equation*}%
As a consequence of our choice of $\left( \hbar ,k\right) $, we conclude
that 
\begin{equation*}
\left\Vert y\left( T_{3}\right) \right\Vert \leq \sqrt{\hbar }\left\Vert
z\right\Vert =\varepsilon \left\Vert z\right\Vert
\end{equation*}%
and 
\begin{equation*}
\left\Vert f\right\Vert _{\omega }\leq \mathcal{C}_{\varepsilon }\left(
T_{3}-T_{2},T_{2}-T_{1}\right) \left\Vert z\right\Vert \text{.}
\end{equation*}%
From these, we see that the above $f$ satisfies two properties in $(i)$ of
Theorem \ref{theorem3.1}.

It remains to prove the above claim. To this end, we write the inequality in
Theorem \ref{theorem2.1} $(iii)$ as follows. 
\begin{equation*}
\left\Vert e^{tA}\Phi \right\Vert ^{2}\leq \left( \left\Vert \Phi
\right\Vert ^{2}\right) ^{\theta }\left( e^{\frac{2C_{3}}{1-\theta }\left( 1+%
\frac{1}{\theta t}+t\left\Vert V\right\Vert _{\infty }+\left\Vert
V\right\Vert _{\infty }^{2/3}\right) }\left\Vert 1_{\omega }^{\ast
}e^{tA}\Phi \right\Vert _{\omega }^{2}\right) ^{1-\theta }\text{ .}
\end{equation*}%
Since $\left\Vert e^{LA}\Phi \right\Vert \leq e^{\left( L-t\right)
\left\Vert V\right\Vert _{\infty }}\left\Vert e^{tA}\Phi \right\Vert $ for $%
L\geq t$, it holds that 
\begin{equation*}
\left\Vert e^{LA}\Phi \right\Vert ^{2}\leq \left( \left\Vert \Phi
\right\Vert ^{2}\right) ^{\theta }\left( e^{\frac{2}{1-\theta }\left(
L-t\right) \left\Vert V\right\Vert _{\infty }}e^{\frac{2C_{3}}{1-\theta }%
\left( 1+\frac{1}{\theta t}+t\left\Vert V\right\Vert _{\infty }+\left\Vert
V\right\Vert _{\infty }^{2/3}\right) }\left\Vert 1_{\omega }^{\ast
}e^{tA}\Phi \right\Vert _{\omega }^{2}\right) ^{1-\theta }\text{ .}
\end{equation*}%
Following the same technique as that used in the proof of Theorem \ref%
{theorem2.1} (step 5), using the Young inequality, we have that for any $%
\varepsilon >0$ and any $\Phi \in L^{2}\left( \Omega \right) $, 
\begin{equation*}
\left\Vert e^{LA}\Phi \right\Vert ^{2}\leq \text{exp}\left( 2\Lambda
_{L}+2\Upsilon +\beta \left( \text{ln}\frac{1}{\varepsilon }+2\Lambda
_{L}\right) +\frac{1}{\beta }2\Upsilon \right) \left\Vert 1_{\omega }^{\ast
}e^{tA}\Phi \right\Vert _{\omega }^{2}+\varepsilon \left\Vert \Phi
\right\Vert ^{2}\text{ ,}
\end{equation*}%
that is equivalent to%
\begin{equation*}
\left\Vert e^{LA}\Phi \right\Vert ^{2}\leq \text{exp}\left( 2\left( \Lambda
_{L}+\Upsilon +\beta \left( \text{ln}\frac{1}{\varepsilon }+\Lambda
_{L}\right) +\frac{1}{\beta }\Upsilon \right) \right) \left\Vert 1_{\omega
}^{\ast }e^{tA}\Phi \right\Vert _{\omega }^{2}+\varepsilon ^{2}\left\Vert
\Phi \right\Vert ^{2}\text{ ,}
\end{equation*}%
with $\Lambda _{L}=\left( L-t\right) \left\Vert V\right\Vert _{\infty
}+C_{3}\left( 1+\frac{1}{t}+t\left\Vert V\right\Vert _{\infty }+\left\Vert
V\right\Vert _{\infty }^{2/3}\right) $ and $\Upsilon =\frac{C_{3}}{t}$.
Next, we choose $\beta =\sqrt{\frac{\Upsilon }{\text{ln}^{+}\frac{1}{%
\varepsilon }+\Lambda _{L}}}$ to get that 
\begin{equation*}
\left\Vert e^{LA}\Phi \right\Vert ^{2}\leq \text{exp}\left( 8\Lambda _{L}+4%
\sqrt{\Upsilon \text{ln}^{+}\frac{1}{\varepsilon }}\right) \left\Vert
1_{\omega }^{\ast }e^{tA}\Phi \right\Vert _{\omega }^{2}+\varepsilon
^{2}\left\Vert \Phi \right\Vert ^{2}\text{ ,}
\end{equation*}%
for any $\varepsilon >0$, $L\geq t>0$ and $\Phi \in L^{2}\left( \Omega
\right) $. Setting $c=4C_{3}$ (here $c>0$ is the same constant than in $(v)$
of Theorem 2.1) and applying the above, with the choices $L=T_{3}-T_{1}$, $%
t=T_{3}-T_{2}$, and $\Phi =w$, give the desired claim.

This completes the proof of the conclusion $\left( i\right) $ in Theorem \ref%
{theorem3.1}.

\bigskip

\bigskip

\bigskip

\subsection{Uniqueness of minimal norm impulse control and its construction}

\bigskip

The aim of this subsection is to study a minimal norm problem, which is
indeed given in $\left( ii\right) $ of Theorem \ref{theorem3.1}. We will
present some properties on this problem (see Theorem~\ref{theorem3.4}), and
give the proof of $\left( ii\right) $ of Theorem \ref{theorem3.1} (see $%
\left( a\right) $ of Remark~\ref{remark3.5}).

Arbitrarily fix $z\in L^{2}\left( \Omega \right) \left\backslash
\{0\}\right. $ and $\varepsilon >0$. Recall the following impulse controlled
equation over $\left[ T_{1},T_{3}\right] $:%
\begin{equation}
\left\{ 
\begin{array}{ll}
y^{\prime }\left( t\right) -Ay\left( t\right) =0\text{ ,} & t\in \left(
T_{1},T_{3}\right) \backslash \left\{ T_{2}\right\} \ \text{,} \\ 
y\left( T_{1}\right) =z\text{ ,} &  \\ 
y\left( T_{2}\right) =y\left( T_{2-}\right) +1_{\omega }f\text{ .} & 
\end{array}%
\right.  \tag{3.2.1}  \label{3.2.1}
\end{equation}%
In this subsection, we discuss the following minimal norm impulse control
problem $(\mathcal{P})$:%
\begin{equation}
\mathcal{N}_{z}:=\text{inf}\left\{ \left\Vert f\right\Vert _{\omega };f\in
L^{2}\left( \omega \right) \text{ and }\left\Vert y\left( T_{3}\right)
\right\Vert \leq \varepsilon \left\Vert z\right\Vert \right\} \text{ .} 
\tag{3.2.2}  \label{3.2.2}
\end{equation}%
This problem is to ask for a control which has the minimal norm among all
controls (in $L^{2}\left( \omega \right) $) driving solutions of equation (%
\ref{3.2.1}) from the initial state $z$ to the closed ball in $L^{2}(\Omega
) $, centered at the origin and of radius $\varepsilon \left\Vert
z\right\Vert $, at the ending time $T_{3}$. In this problem, $\mathcal{N}%
_{z} $ is called the minimal norm, while $f^{\ast }\in L^{2}\left( \omega
\right) $ is called a minimal norm control, if the solution $y$ of (\ref%
{3.2.1}) with $f=f^{\ast }$ satisfies 
\begin{equation*}
\left\Vert y\left( T_{3}\right) \right\Vert \leq \varepsilon \left\Vert
z\right\Vert \text{ and }\left\Vert f^{\ast }\right\Vert _{\omega }=\mathcal{%
N}_{z}\text{ .}
\end{equation*}

\bigskip

The main result of this subsection is as:

\bigskip

\begin{thm}
\label{theorem3.4} The following conclusions are true:

\begin{description}
\item[$\left( i\right) $] The problem $(\mathcal{P})$ has a unique minimal
norm control.

\item[$\left( ii\right) $] The minimal norm control $f^{\ast }$\ to $(%
\mathcal{P})$ satisfies that 
\begin{equation*}
f^{\ast }=0\text{ if and only if the solution }y^{0}\text{ of (\ref{3.2.1})
with }f=0\text{ satisfies }\left\Vert y^{0}\left( T_{3}\right) \right\Vert
\leq \varepsilon \left\Vert z\right\Vert \text{ .}
\end{equation*}

\item[$\left( iii\right) $] The minimal norm control $f^{\ast }$\ to $(%
\mathcal{P})$ is given by 
\begin{equation*}
f^{\ast }=1_{\omega }^{\ast }e^{\left( T_{3}-T_{2}\right) }w
\end{equation*}%
where $w$ is the unique minimizer to $J:L^{2}\left( \Omega \right)
\rightarrow \mathbb{R}$ defined by 
\begin{equation*}
J\left( \Phi \right) :=\frac{1}{2}\left\Vert 1_{\omega }^{\ast }\,e^{\left(
T_{3}-T_{2}\right) A}\Phi \right\Vert _{\omega }^{2}+\left\langle
z,e^{\left( T_{3}-T_{1}\right) A}\Phi \right\rangle +\varepsilon \left\Vert
z\right\Vert \left\Vert \Phi \right\Vert \text{ .}
\end{equation*}
\end{description}
\end{thm}

\bigskip

\begin{rem}
\label{remark3.5} $\left( a\right) $ The conclusion $(i)$ of Theorem \ref%
{theorem3.4} clearly gives the uniqueness in $(ii)$ of Theorem \ref%
{theorem3.1}. This, along with $(i)$ of Theorem \ref{theorem3.1}, shows the
conclusion $(ii)$ in Theorem \ref{theorem3.1}. $\left( b\right) $ The
construction of the control in $(iii)$ of Theorem \ref{theorem3.4} is
inspired by a standard duality strategy used in \cite{FPZ} for the
distributed controlled heat equations with a control in $L^{2}\left( \omega
\times \left( T_{1},T_{3}\right) \right) $.
\end{rem}

\bigskip

Proof of Theorem \ref{theorem3.4} .-

\noindent \textit{Proof of }$(i)$\textit{: }Write%
\begin{equation*}
\mathcal{F}_{ad}:=\left\{ f\in L^{2}\left( \omega \right) ;\left\Vert
y\left( T_{3}\right) \right\Vert \leq \varepsilon \left\Vert z\right\Vert
\right\} \text{ .}
\end{equation*}%
By the previous subsection, we see that $\mathcal{F}_{ad}\neq \emptyset $.
Meanwhile, one can easily check that $\mathcal{F}_{ad}$ is weakly closed in $%
L^{2}\left( \omega \right) $. From these, it follows that $(\mathcal{P})$
has a minimal norm control.

Suppose that $f_{1}$ and $f_{2}$ are two minimal norm controls to $(\mathcal{%
P})$. Then we have that 
\begin{equation*}
0\leq \Vert f_{1}\Vert _{\omega }=\Vert f_{2}\Vert _{\omega }=\mathcal{N}%
_{z}<\infty \text{ ,}
\end{equation*}%
where $\mathcal{N}_{z}$ is given by (\ref{3.2.2}). Meanwhile, one can easily
check that $(f_{1}+f_{2})/2$ is also a minimal norm control to $(\mathcal{P}%
) $. This, along with the Parallelogram Law, yields that 
\begin{equation*}
\left( \mathcal{N}_{z}\right) ^{2}=\Vert (f_{1}+f_{2})/2\Vert _{\omega }^{2}=%
\frac{1}{2}\left( \Vert f_{1}\Vert _{\omega }^{2}+\Vert f_{2}\Vert _{\omega
}^{2}\right) -\Vert (f_{1}-f_{2})/2\Vert _{\omega }^{2}\text{ .}
\end{equation*}%
From the two above identities on $\mathcal{N}_{z}$, we find that $%
f_{1}=f_{2} $. Thus, the minimal norm control to $(\mathcal{P})$ is unique.
This ends the proof.

\noindent \textit{Proof of }$(ii)$\textit{: }The second conclusion in
Theorem \ref{theorem3.4} follows from the definition of Problem $(\mathcal{P}%
)$ (see (\ref{3.2.2})) at once.

\noindent \textit{Proof of }$(iii)$\textit{:} To prove the last conclusion
in Theorem \ref{theorem3.4}, we need the next Lemma \ref{lemma3.6} whose
proof will be given at the end of the proof of Theorem \ref{theorem3.4}.

\bigskip

\begin{lem}
\label{lemma3.6} The functional $J$ in Theorem \ref{theorem3.4} has the
following properties:

\begin{description}
\item[$\left( a\right) $] It satisfies that 
\begin{equation}
\underset{q\rightarrow \infty }{\normalfont{\text{lim}}}\underset{\left\Vert
\Phi \right\Vert =q}{\normalfont{\text{inf}}}\frac{J\left( \Phi \right) }{%
\left\Vert \Phi \right\Vert }\geq \varepsilon \left\Vert z\right\Vert \text{
.}  \tag{3.2.3}  \label{3.2.3}
\end{equation}

\item[$\left( b\right) $] It has a unique minimizer over $L^{2}(\Omega )$ .

\item[$\left( c\right) $] Write $w$ for its minimizer. Then%
\begin{equation*}
w=0\text{ if and only if the solution }y^{0}\text{ of (\ref{3.2.1}) with }f=0%
\text{ satisfies }\left\Vert y^{0}\left( T_{3}\right) \right\Vert \leq
\varepsilon \left\Vert z\right\Vert \text{ .}
\end{equation*}
\end{description}
\end{lem}

\bigskip

We now show the third conclusion of Theorem \ref{theorem3.4}. Notice that
when $\left\Vert y^{0}\left( T_{3}\right) \right\Vert \leq \varepsilon
\left\Vert z\right\Vert $ where $y^{0}$ is the solution of (\ref{3.2.1})
with $f=0$, it follows respectively from the second conclusion of Theorem %
\ref{theorem3.4} and the conclusion $(c)$ of Lemma \ref{lemma3.6} that $%
f^{\ast }=0$ and $w=0$. Hence, the third conclusion of Theorem 3.4 is true
in this particular case.

We now consider the case where 
\begin{equation}
\left\Vert y^{0}\left( T_{3}\right) \right\Vert >\varepsilon \left\Vert
z\right\Vert  \tag{3.2.4}  \label{3.2.4}
\end{equation}%
where $y^{0}$ is the solution of (\ref{3.2.1}) with $f=0$. Let $w\in
L^{2}(\Omega )$ be the minimizer of the functional $J$ (see $(b)$ of Lemma %
\ref{lemma3.6}). Write 
\begin{equation}
\widehat{f}=1_{\omega }^{\ast }e^{\left( T_{3}-T_{2}\right) A}w\text{ .} 
\tag{3.2.5}  \label{3.2.5}
\end{equation}%
We first claim that 
\begin{equation}
\widehat{f}\in \mathcal{F}_{ad}:=\left\{ f\in L^{2}\left( \omega \right)
;\left\Vert y\left( T_{3}\right) \right\Vert \leq \varepsilon \left\Vert
z\right\Vert \right\} \text{ .}  \tag{3.2.6}  \label{3.2.6}
\end{equation}%
In fact, by (\ref{3.2.4}) and $(c)$ of Lemma \ref{lemma3.6}, we find that $%
w\neq 0$. Then the Euler-Lagrange equation associated to $w$ reads:%
\begin{equation}
e^{\left( T_{3}-T_{2}\right) A}\chi _{\omega }e^{\left( T_{3}-T_{2}\right)
A}w+e^{\left( T_{3}-T_{1}\right) A}z+\varepsilon \left\Vert z\right\Vert 
\frac{w}{\left\Vert w\right\Vert }=0\text{ .}  \tag{3.2.7}  \label{3.2.7}
\end{equation}%
Meanwhile, since $1_{\omega }1_{\omega }^{\ast }=\chi _{\omega }$, it
follows from (\ref{3.2.5}) that the solution $y$ of (\ref{3.2.1}) with $f=%
\widehat{f}$ satisfies 
\begin{equation*}
y\left( T_{3}\right) =e^{\left( T_{3}-T_{1}\right) A}z+e^{\left(
T_{3}-T_{2}\right) A}1_{\omega }\widehat{f}=e^{\left( T_{3}-T_{2}\right)
A}\chi _{\omega }e^{\left( T_{3}-T_{2}\right) A}w+e^{\left(
T_{3}-T_{1}\right) A}z\text{ .}
\end{equation*}%
This, together with (\ref{3.2.7}), indicates that 
\begin{equation*}
y\left( T_{3}\right) =-\varepsilon \left\Vert z\right\Vert \frac{w}{%
\left\Vert w\right\Vert }\text{ ,}
\end{equation*}%
from which, (\ref{3.2.6}) follows at once.

We next claim that 
\begin{equation}
\left\Vert \widehat{f}\right\Vert _{\omega }\leq \left\Vert f\right\Vert
_{\omega }\quad \text{for all }f\in \mathcal{F}_{ad}\text{ .}  \tag{3.2.8}
\label{3.2.8}
\end{equation}%
To this end, we arbitrarily fix an $f\in \mathcal{F}_{ad}$. Then we have
that the solution $y$ of (\ref{3.2.1}) satisfies 
\begin{equation}
\left\Vert y\left( T_{3}\right) \right\Vert \leq \varepsilon \left\Vert
z\right\Vert \text{ .}  \tag{3.2.9}  \label{3.2.9}
\end{equation}%
Since $\widehat{f}:=1_{\omega }^{\ast }e^{\left( T_{3}-T_{2}\right) A}w$, 
\begin{equation}
\begin{array}{ll}
\left\Vert \widehat{f}\right\Vert _{\omega }^{2} & =\left\Vert 1_{\omega
}^{\ast }e^{\left( T_{3}-T_{2}\right) A}w\right\Vert _{\omega }^{2} \\ 
& =\left\Vert 1_{\omega }^{\ast }e^{\left( T_{3}-T_{2}\right) A}w\right\Vert
_{\omega }^{2}+\varepsilon \left\Vert z\right\Vert \left\Vert w\right\Vert
-\varepsilon \left\Vert z\right\Vert \left\Vert w\right\Vert \\ 
& =-\varepsilon \left\Vert z\right\Vert \left\Vert w\right\Vert
+\left\langle e^{\left( T_{3}-T_{2}\right) A}\chi _{\omega }e^{\left(
T_{3}-T_{2}\right) A}w+\varepsilon \left\Vert z\right\Vert \frac{w}{%
\left\Vert w\right\Vert },w\right\rangle \\ 
& =-\varepsilon \left\Vert z\right\Vert \left\Vert w\right\Vert
-\left\langle e^{\left( T_{3}-T_{1}\right) A}z,w\right\rangle \\ 
& \leq \left\langle y\left( T_{3}\right) ,w\right\rangle -\left\langle
e^{\left( T_{3}-T_{1}\right) A}z,w\right\rangle =\left\langle e^{\left(
T_{3}-T_{2}\right) A}1_{\omega }f,w\right\rangle \\ 
& =\left\langle 1_{\omega }f,\chi _{\omega }e^{\left( T_{3}-T_{2}\right)
A}w\right\rangle =\left\langle 1_{\omega }f,\,1_{\omega }\widehat{f}%
\right\rangle \leq \frac{1}{2}\left\Vert f\right\Vert _{\omega }^{2}+\frac{1%
}{2}\left\Vert \widehat{f}\right\Vert _{\omega }^{2}\text{ .}%
\end{array}
\tag{3.2.10}  \label{3.2.10}
\end{equation}
Notice that we used $1_{\omega }1_{\omega }^{\ast }=\chi _{\omega }$, $\chi
_{\omega }1_{\omega }=1_{\omega }$ and Cauchy-Schwarz inequality in the last
line in (\ref{3.2.10}); the equality (\ref{3.2.7}) is applied in the fourth
equality of (\ref{3.2.10}); the inequality (\ref{3.2.9}) and Cauchy-Schwarz,
as well as the formula $y\left( T_{3}\right) =e^{\left( T_{3}-T_{1}\right)
A}z+e^{\left( T_{3}-T_{2}\right) A}1_{\omega }f$, are used in the fifth line
of (\ref{3.2.10}). Now, (\ref{3.2.10}) clearly leads to (\ref{3.2.8}).

From (\ref{3.2.6}) and (\ref{3.2.8}), we find that $\widehat{f}$ is a
minimal norm control to $(\mathcal{P})$. Since the minimal norm control of $(%
\mathcal{P})$ is unique, we have that $\widehat{f}=f^{\ast }$. So the third
conclusion of Theorem 3.4 is true.

\bigskip

Finally, we are ready to prove Lemma \ref{lemma3.6}.

\bigskip

Proof of Lemma \ref{lemma3.6} .- We will prove conclusions $(a),\left(
b\right) ,(c)$ one by one.

\noindent \textit{Proof of }$(a)$\textit{:} By contradiction, suppose that (%
\ref{3.2.3}) was not true. Then there would be an $\sigma \in (0,\varepsilon
)$ and a sequence $\{\Phi _{n}\}_{n=1}^{\infty }$ in $L^{2}(\Omega )$ so
that 
\begin{equation}
\underset{n\rightarrow \infty }{\text{lim}}\left\Vert \Phi _{n}\right\Vert
=\infty \text{ }  \tag{3.2.11}  \label{3.2.11}
\end{equation}%
and%
\begin{equation}
\frac{J\left( \Phi _{n}\right) }{\left\Vert \Phi _{n}\right\Vert }\leq
\left( \varepsilon -\sigma \right) \left\Vert z\right\Vert \text{\quad for
all }n\in \mathbb{N}\text{ .}  \tag{3.2.12}  \label{3.2.12}
\end{equation}%
From (\ref{3.2.11}), we can assume, without loss of generality, that $\Phi
_{n}\neq 0$ for all $n$. Thus we can set 
\begin{equation}
\varphi _{n}=\frac{\Phi _{n}}{\Vert \Phi _{n}\Vert }\quad \text{for all }%
n\in \mathbb{N}\text{ .}  \tag{3.2.13}  \label{3.2.13}
\end{equation}%
From (\ref{3.2.13}), we see that $\left\{ e^{\left( T_{3}-T_{1}\right)
A}\varphi _{n}\right\} _{n=1}^{\infty }$ is bounded in $L^{2}(\Omega )$.
Then, from the definition of $J$ in Theorem \ref{theorem3.4} , (\ref{3.2.13}%
), (\ref{3.2.11}) and (\ref{3.2.12}), we find that 
\begin{equation}
\begin{array}{ll}
& \quad \underset{n\rightarrow \infty }{\overline{\text{lim}}}\displaystyle%
\frac{1}{2}\left\Vert 1_{\omega }^{\ast }\,e^{\left( T_{3}-T_{2}\right)
A}\varphi _{n}\right\Vert _{\omega }^{2} \\ 
& =\underset{n\rightarrow \infty }{\overline{\text{lim}}}\displaystyle\frac{1%
}{\Vert \Phi _{n}\Vert }\left[ \displaystyle\frac{J\left( \Phi _{n}\right) }{%
\left\Vert \Phi _{n}\right\Vert }-\left\langle z,e^{\left(
T_{3}-T_{1}\right) A}\varphi _{n}\right\rangle -\varepsilon \left\Vert
z\right\Vert \right] \\ 
& \leq \underset{n\rightarrow \infty }{\overline{\text{lim}}}\displaystyle%
\frac{-\sigma \left\Vert z\right\Vert }{\Vert \Phi _{n}\Vert }+\underset{%
n\rightarrow \infty }{\overline{\text{lim}}}\displaystyle\frac{-\left\langle
z,e^{\left( T_{3}-T_{1}\right) A}\varphi _{n}\right\rangle }{\Vert \Phi
_{n}\Vert }=0\text{ .}%
\end{array}
\tag{3.2.14}  \label{3.2.14}
\end{equation}%
Meanwhile, by (\ref{3.2.13}), there is a subsequence of $\left\{ \varphi
_{n}\right\} $, denoted in the same manner, so that 
\begin{equation*}
\varphi _{n}\rightarrow \varphi \quad \text{weakly in }L^{2}\left( \Omega
\right) \text{ ,}
\end{equation*}%
for some $\varphi \in L^{2}\left( \Omega \right) $. Since the semigroup $%
\left\{ e^{tA}\right\} _{t\geq 0}$ is compact, the above convergence leads to%
\begin{equation}
e^{\left( T_{3}-T_{2}\right) A}\varphi _{n}\rightarrow e^{\left(
T_{3}-T_{2}\right) A}\varphi \quad \text{strongly in }L^{2}\left( \Omega
\right) \text{ }  \tag{3.2.15}  \label{3.2.15}
\end{equation}%
and 
\begin{equation}
1_{\omega }^{\ast }e^{\left( T_{3}-T_{2}\right) A}\varphi _{n}\rightarrow
1_{\omega }^{\ast }e^{\left( T_{3}-T_{2}\right) A}\varphi \quad \text{%
strongly in }L^{2}\left( \omega \right) \text{ .}  \tag{3.2.16}
\label{3.2.16}
\end{equation}%
From (\ref{3.2.14}) and the convergence in (\ref{3.2.16}), we find that 
\begin{equation*}
1_{\omega }^{\ast }e^{\left( T_{3}-T_{2}\right) A}\varphi =0\text{ ,}
\end{equation*}%
which, along with the unique continuation property of heat equations (see $%
(1)$ of Remark \ref{remark2.2}) yields that $\varphi =0$. Then from the
definition of $J$ in Theorem \ref{theorem3.4} and the convergence in (\ref%
{3.2.15}), we see that 
\begin{equation*}
\begin{array}{ll}
\underset{n\rightarrow \infty }{\underline{\text{lim}}}\displaystyle\frac{%
J\left( \Phi _{n}\right) }{\left\Vert \Phi _{n}\right\Vert } & \geq \underset%
{n\rightarrow \infty }{\underline{\text{lim}}}\left[ \left\langle
z,e^{\left( T_{3}-T_{1}\right) A}\varphi _{n}\right\rangle +\varepsilon
\left\Vert z\right\Vert \right] \\ 
& =\left\langle z,e^{\left( T_{3}-T_{1}\right) A}\varphi \right\rangle
+\varepsilon \left\Vert z\right\Vert =\varepsilon \left\Vert z\right\Vert 
\text{ .}%
\end{array}%
\end{equation*}%
This, along with (\ref{3.2.12}), leads to a contradiction. Therefore, (\ref%
{3.2.3}) is true.

\noindent \textit{Proof of }$(b)$\textit{:} From (\ref{3.2.3}), we see that
the functional $J$ is coercive on $L^{2}(\Omega )$. Further $J$ is
continuous and convex on $L^{2}(\Omega )$. Thus, it has a minimizer on $%
L^{2}(\Omega )$.

Next, we show the uniqueness of the minimizer. It suffices to prove that the
functional $J$ is strictly convex. For this purpose, we arbitrarily fix $%
\Phi _{1},\Phi _{2}\in L^{2}\left( \Omega \right) \left\backslash
\{0\}\right. $, with $\Phi _{1}\neq \Phi _{2}$. There are only three
possibilities: $(a)$ $\Phi _{1}\neq \mu \Phi _{2}$ for any $\mu \in \mathbb{R%
}$; $(b)$ $\Phi _{1}=-\mu _{0}\Phi _{2}$ for some $\mu _{0}>0$; $(c)$ $\Phi
_{1}=\mu _{0}\Phi _{2}$ for some $\mu _{0}>0$. In the cases $(a)$ and $(b)$,
one can easily check that 
\begin{equation}
\Vert \lambda \Phi _{1}+(1-\lambda )\Phi _{2}\Vert <\lambda \Vert \Phi
_{1}\Vert +(1-\lambda )\Vert \Phi _{2}\Vert \quad \text{for all }\lambda \in
(0,1)\text{ .}  \tag{3.2.17}  \label{3.2.17}
\end{equation}%
In the case $(c)$, we let 
\begin{equation*}
H(\lambda )=J\left( \lambda \Phi _{2}\right) \text{ ,}\quad \lambda >0\text{
.}
\end{equation*}%
Since $\Phi _{2}\neq 0$ in $L^{2}(\Omega )$, it follows by the unique
continuation property of heat equations (see $(1)$ of Remark \ref{remark2.2}%
) that $\left\Vert 1_{\omega }^{\ast }e^{\left( T_{3}-T_{2}\right) A}\Phi
_{2}\right\Vert _{\omega }\neq 0$. Thus, $H$ is a quadratic function with a
positive leading coefficient. Hence, $H$ is strictly convex. This, along
with (\ref{3.2.17}), yields the strict convexity of $J$.

\noindent \textit{Proof of }$(c)$\textit{:} Let $y^{0}$ be the solution of (%
\ref{3.2.1}) with $f=0$. We first show that 
\begin{equation}
\left\Vert y^{0}\left( T_{3}\right) \right\Vert \leq \varepsilon \left\Vert
z\right\Vert \Rightarrow w=0\text{ .}  \tag{3.2.18}  \label{3.2.18}
\end{equation}%
In fact, by multiplying by $e^{\left( T_{3}-t\right) A}\Phi $ the system
solved by $y^{0}$, we see that 
\begin{equation*}
\left\langle z,e^{\left( T_{3}-T_{1}\right) A}\Phi \right\rangle
=\left\langle y^{0}\left( T_{3}\right) ,\Phi \right\rangle \quad \text{for
all}\;z,\Phi \in L^{2}(\Omega )\text{ .}
\end{equation*}%
This, along with the definition of $J$ in Theorem \ref{theorem3.4} and the
inequality on the left hand side of (\ref{3.2.18}), yields that for all$%
\;\Phi \in L^{2}(\Omega )$ 
\begin{equation*}
J\left( \Phi \right) \geq \left\langle y^{0}\left( T_{3}\right) ,\Phi
\right\rangle +\varepsilon \left\Vert z\right\Vert \left\Vert \Phi
\right\Vert \geq 0=J\left( 0\right) \text{ .}
\end{equation*}%
This implies the equality on the right hand side of (\ref{3.2.18}). Hence, (%
\ref{3.2.18}) is true.

We next show that 
\begin{equation}
w=0\Rightarrow \left\Vert y^{0}\left( T_{3}\right) \right\Vert \leq
\varepsilon \left\Vert z\right\Vert \text{ .}  \tag{3.2.19}  \label{3.2.19}
\end{equation}%
By contradiction, suppose that (\ref{3.2.19}) were not true. Then we would
have that 
\begin{equation}
\left\Vert y^{0}\left( T_{3}\right) \right\Vert >\varepsilon \left\Vert
z\right\Vert \text{ and }w=0\text{ .}  \tag{3.2.20}  \label{3.2.20}
\end{equation}%
Set $\psi :=-y^{0}\left( T_{3}\right) $, which clearly belongs to $%
L^{2}\left( \Omega \right) \left\backslash \{0\}\right. $. Then we have that 
\begin{equation*}
\left\langle z,e^{\left( T_{3}-T_{1}\right) A}\psi \right\rangle
=\left\langle y^{0}\left( T_{3}\right) ,\psi \right\rangle =-\left\Vert
y^{0}\left( T_{3}\right) \right\Vert \left\Vert \psi \right\Vert \text{ .}
\end{equation*}%
This, along with the first inequality in (\ref{3.2.20}), yields that 
\begin{equation*}
\left\langle z,e^{\left( T_{3}-T_{1}\right) A}\psi \right\rangle
+\varepsilon \left\Vert z\right\Vert \left\Vert \psi \right\Vert <0\text{ .}
\end{equation*}%
Thus, there is an $\sigma >0$ so that 
\begin{equation*}
J\left( \sigma \psi \right) =\sigma ^{2}\frac{1}{2}\left\Vert 1_{\omega
}^{\ast }e^{\left( T_{3}-T_{2}\right) A}\psi \right\Vert _{\omega
}^{2}+\sigma \left( \left\langle z,e^{\left( T_{3}-T_{1}\right) A}\psi
\right\rangle +\varepsilon \left\Vert z\right\Vert \left\Vert \psi
\right\Vert \right) <0\text{ .}
\end{equation*}%
This, along with the second equation in (\ref{3.2.20}), indicates that 
\begin{equation*}
0=J\left( 0\right) =\underset{\Phi \in L^{2}(\Omega )}{\text{min}}J\left(
\Phi \right) <0\text{ ,}
\end{equation*}%
which leads to a contradiction. So we have proved (\ref{3.2.19}). Finally,
the conclusion $\left( c\right) $ of Lemma \ref{lemma3.6} follows from (\ref%
{3.2.18}) and (\ref{3.2.19}) at once.

This ends the proof of Lemma \ref{lemma3.6} and completes the proof of
Theorem \ref{theorem3.4}.

\bigskip

\bigskip

\subsection{Best connection between Theorem 3.4 and Theorem 2.1}

\bigskip

What we study in this subsection will not have influence on the study of our
stabilization. However, it is independently interesting. Consider the
following problem $(\mathcal{NP})$ (with arbitrarily fixed $\varepsilon >0$.
Recall that $y$ is the solution of (\ref{3.2.1}) associated with the initial
datum $z$ and control $f$): 
\begin{equation}
\mathcal{N}:=\underset{\left\Vert z\right\Vert \leq 1}{\text{sup}}\mathcal{N}%
_{z}=\underset{\left\Vert z\right\Vert \leq 1}{\text{sup}}\text{inf}\left\{
\left\Vert f\right\Vert _{\omega };f\in L^{2}\left( \omega \right) \text{
and }\left\Vert y\left( T_{3}\right) \right\Vert \leq \varepsilon \left\Vert
z\right\Vert \right\} \text{ .}  \tag{3.3.1}  \label{3.3.1}
\end{equation}%
The quantity $\mathcal{N}$ is called the value of the problem $(\mathcal{NP})
$. Next, let $C>0$ and introduce the following property $(\mathcal{Q}_{C})$:
For any $z\in L^{2}\left( \Omega \right) $, there is a control $f\in
L^{2}\left( \omega \right) $ so that 
\begin{equation}
\text{max}\left\{ \frac{1}{C}\left\Vert f\right\Vert _{\omega },\frac{1}{%
\varepsilon }\left\Vert y\left( T_{3}\right) \right\Vert \right\} \leq
\left\Vert z\right\Vert \text{ .}  \tag{3.3.2}  \label{3.3.2}
\end{equation}%
We would like to mention that the property $(\mathcal{Q}_{C})$ may not hold
for some $C>0$ and $\varepsilon >0$. However, we have seen in Theorem \ref%
{theorem3.1} that given $\varepsilon >0$, there is $C=\mathcal{C}%
_{\varepsilon }\left( T_{3}-T_{2},T_{2}-T_{1}\right) >0$ so that the
property $(\mathcal{Q}_{C})$ is true.

\bigskip

The main result of this subsection is as follows: The value $\mathcal{N}$ is
the optimal coefficient $C$ so that%
\begin{equation}
\left\Vert e^{\left( T_{3}-T_{1}\right) A}\Phi \right\Vert \leq C\left\Vert
1_{\omega }^{\ast }e^{\left( T_{3}-T_{2}\right) A}\Phi \right\Vert _{\omega
}+\varepsilon \left\Vert \Phi \right\Vert \quad \text{for any }\Phi \in
L^{2}\left( \Omega \right) \text{ .}  \tag{3.3.3}  \label{3.3.3}
\end{equation}%
Precisely, we have the following result:

\bigskip

\begin{thm}
\label{theorem3.7} Let $\varepsilon >0$. It holds that 
\begin{equation*}
\normalfont{\text{inf}}\left\{ C>0;C\text{ \textit{satisfies (\ref{3.3.3})}}%
\right\} =\mathcal{N}\text{ .}
\end{equation*}
\end{thm}

\bigskip

Further the connections among the problem $(\mathcal{NP})$, the property $(%
\mathcal{Q}_{C})$ and the observation inequalities in Theorem \ref%
{theorem2.1} are presented in the next Theorem~\ref{theorem3.8}, which will
be used in the proof of the above Theorem~\ref{theorem3.7}.

\bigskip

\begin{thm}
\label{theorem3.8} Let $\varepsilon >0$ and $C>0$. The following statements
are equivalent:

\begin{description}
\item[$\left( i\right) $] Let $\mathcal{N}$ be given by (\ref{3.3.1}). Then $%
\mathcal{N}\leq C$.

\item[$\left( ii\right) $] The property $(\mathcal{Q}_{C})$ defined by (\ref%
{3.3.2}) is true.

\item[$\left( iii\right) $] For any $\Phi \in L^{2}\left( \Omega \right) $,
the following estimate holds: 
\begin{equation*}
\left\Vert e^{\left( T_{3}-T_{1}\right) A}\Phi \right\Vert \leq C\left\Vert
1_{\omega }^{\ast }e^{\left( T_{3}-T_{2}\right) A}\Phi \right\Vert _{\omega
}+\varepsilon \left\Vert \Phi \right\Vert \text{ .}
\end{equation*}
\end{description}
\end{thm}

\bigskip

Proof of Theorem \ref{theorem3.8} .- We organize the proof by three steps as
follows:

\noindent\textit{Step 1. To show that }$(i)\Leftrightarrow(ii)$.

We first prove that $(i)\Rightarrow (ii)$. Assume that $(i)$ is true. When $%
z=0$ in $L^{2}(\Omega )$, we find that (\ref{3.3.2}) holds for $f=0$. Thus,
it suffices to show $(ii)$ with an arbitrarily fixed $z\in L^{2}\left(
\Omega \right) \left\backslash \{0\}\right. $. For this purpose, we write $%
\widehat{z}=z/\left\Vert z\right\Vert $. Let $\widehat{f}$ be the solution
to $(\mathcal{P})$ associated to $\widehat{z}$. Then the solution $\widehat{y%
}$ of (\ref{3.2.1}) associated with initial data $\widehat{z}$ and control $%
\widehat{f}$ satisfies that 
\begin{equation*}
\left\Vert \widehat{y}\left( T_{3}\right) \right\Vert \leq \varepsilon
\left\Vert \widehat{z}\right\Vert =\varepsilon \text{ .}
\end{equation*}%
Setting $f=\left\Vert z\right\Vert \widehat{f}$, the solution $y$ of (\ref%
{3.2.1}) have the following property: 
\begin{equation*}
\left\Vert y\left( T_{3}\right) \right\Vert =\left\Vert z\right\Vert
\left\Vert \widehat{y}\left( T_{3}\right) \right\Vert \leq \varepsilon
\left\Vert z\right\Vert \text{ .}
\end{equation*}%
Thus, to show that the above $f$ satisfies (\ref{3.3.2}), we only need to
prove that $\left\Vert f\right\Vert _{\omega }\leq C\left\Vert z\right\Vert $%
. This will be done in what follows: Since $\widehat{f}$ is the solution to $%
(\mathcal{P})$ associated to $\widehat{z}$, we have that $\left\Vert 
\widehat{f}\right\Vert _{\omega }=\mathcal{N}_{\widehat{z}}$. This, along
with (\ref{3.3.1}) and $(i)$ of Theorem \ref{theorem3.8}, yields that 
\begin{equation*}
\left\Vert f\right\Vert _{\omega }=\left\Vert z\right\Vert \left\Vert 
\widehat{f}\right\Vert _{\omega }=\left\Vert z\right\Vert \mathcal{N}_{%
\widehat{z}}\leq \left\Vert z\right\Vert \mathcal{N}\leq C\left\Vert
z\right\Vert \text{ .}
\end{equation*}%
Hence, $(ii)$ is true.

We next show that $(ii)\Rightarrow (i)$. Assume that $(ii)$ is true. By
contradiction, suppose that $(i)$ were false. Then there would be $z\neq 0$
with $\left\Vert z\right\Vert \leq 1$ so that $\mathcal{N}_{z}>C$. Let $%
\widehat{z}=z\left/ \left\Vert z\right\Vert \right. $. Then we have that 
\begin{equation*}
\mathcal{N}_{\widehat{z}}=\frac{1}{\Vert z\Vert }\mathcal{N}_{z}\geq 
\mathcal{N}_{z}>C\text{ .}
\end{equation*}%
Therefore, we see that there is no $f\in L^{2}(\Omega )$ so that the
solution $y$ of (\ref{3.2.1}), associated with the initial datum $\widehat{z}
$ and the control $f$, has the property: 
\begin{equation*}
\left\Vert y\left( T_{3}\right) \right\Vert \leq \varepsilon \left\Vert 
\widehat{z}\right\Vert \text{ and }\left\Vert f\right\Vert _{\omega }\leq
C=C\left\Vert \widehat{z}\right\Vert \text{ .}
\end{equation*}%
This contradicts $(ii)$. Hence, $(i)$ stands.

\noindent\textit{Step 2. To show that }$(ii)\Rightarrow(iii)$.

Suppose that $(ii)$ holds. Then, given $z\in L^{2}(\Omega )$, there is $f\in
L^{2}(\Omega )$ so that (\ref{3.3.2}) holds. Meanwhile, by multiplying by $%
e^{\left( T_{3}-t\right) A}\Phi $ the system solved by $y$, one gets 
\begin{equation*}
\left\langle y\left( T_{3}\right) ,\Phi \right\rangle -\left\langle
z,e^{\left( T_{3}-T_{1}\right) A}\Phi \right\rangle =\left\langle 1_{\omega
}f,e^{\left( T_{3}-T_{2}\right) A}\Phi \right\rangle \quad \forall \Phi \in
L^{2}\left( \Omega \right) \text{ .}
\end{equation*}%
This, along with the inequality (\ref{3.3.2}), yields that for each $\Phi
\in L^{2}(\Omega )$,%
\begin{equation*}
\begin{array}{ll}
& \quad \left\Vert e^{\left( T_{3}-T_{1}\right) A}\Phi \right\Vert =\underset%
{\left\Vert z\right\Vert \leq 1}{\text{sup}}\left\langle e^{\left(
T_{3}-T_{1}\right) A}\Phi ,z\right\rangle \\ 
& =\underset{\left\Vert z\right\Vert \leq 1}{\text{sup}}\left[ \left\langle
y\left( T_{3}\right) ,\Phi \right\rangle -\left\langle 1_{\omega
}f,e^{\left( T_{3}-T_{2}\right) A}\Phi \right\rangle \right] \\ 
& \leq \underset{\left\Vert z\right\Vert \leq 1}{\text{sup}}\left[
\left\Vert y\left( T_{3}\right) \right\Vert \left\Vert \Phi \right\Vert
+\left\Vert f\right\Vert _{\omega }\left\Vert 1_{\omega }^{\ast }e^{\left(
T_{3}-T_{2}\right) A}\Phi \right\Vert _{\omega }\right] \\ 
& \leq \underset{\left\Vert z\right\Vert \leq 1}{\text{sup}}\left[ \left(
\varepsilon \left\Vert \Phi \right\Vert +C\left\Vert 1_{\omega }^{\ast
}e^{\left( T_{3}-T_{2}\right) A}\Phi \right\Vert _{\omega }\right)
\left\Vert z\right\Vert \right] \\ 
& =C\left\Vert 1_{\omega }^{\ast }e^{\left( T_{3}-T_{2}\right) A}\Phi
\right\Vert _{\omega }+\varepsilon \left\Vert \Phi \right\Vert \text{ ,}%
\end{array}%
\end{equation*}%
which leads to the desired observation estimate. Hence, $(iii)$ is true.

\noindent\textit{Step 3. To show that }$(iii)\Rightarrow(ii)$.

Suppose that $(iii)$ is true. Arbitrarily fix $z\in L^{2}\left( \Omega
\right) $. Denote by $y^{0}$ the solution of (\ref{3.2.1}) with $f=0$. In
the case that $\left\Vert y^{0}\left( T_{3}\right) \right\Vert \leq
\varepsilon \left\Vert z\right\Vert $, (\ref{3.3.2}) holds for $f=0$. Thus,
we only need to consider the case that 
\begin{equation}
\left\Vert y^{0}\left( T_{3}\right) \right\Vert >\varepsilon \left\Vert
z\right\Vert \text{ .}  \tag{3.3.4}  \label{3.3.4}
\end{equation}%
In this case, we let $f:=1_{\omega }^{\ast }e^{\left( T_{3}-T_{2}\right) A}w$%
, where $w$ is the unique minimizer of the functional $J$, which is given in
Theorem \ref{theorem3.4}. Then according to $(iii)$ of Theorem \ref%
{theorem3.4}, $f$ is the minimal norm control to $(\mathcal{P})$. By Lemma %
\ref{lemma3.6} and (\ref{3.3.4}), we see that $w\neq 0$. Then using the
Euler-Lagrange equation (\ref{3.2.7}) and noticing that $\chi _{\omega
}=1_{\omega }1_{\omega }^{\ast }$, we find that%
\begin{equation*}
\begin{array}{ll}
& \quad \left\langle z,e^{\left( T_{3}-T_{1}\right) A}w\right\rangle
=\left\langle e^{\left( T_{3}-T_{1}\right) A}z,w\right\rangle  \\ 
& =-\left\langle e^{\left( T_{3}-T_{2}\right) A}\chi _{\omega }e^{\left(
T_{3}-T_{2}\right) A}w+\varepsilon \left\Vert z\right\Vert \frac{w}{%
\left\Vert w\right\Vert },w\right\rangle  \\ 
& =-\left\Vert 1_{\omega }^{\ast }e^{\left( T_{3}-T_{2}\right)
A}w\right\Vert _{\omega }^{2}-\varepsilon \left\Vert z\right\Vert \left\Vert
w\right\Vert \text{ .}%
\end{array}%
\end{equation*}%
Since $f=1_{\omega }^{\ast }e^{\left( T_{3}-T_{2}\right) A}w$, the above
equality, along with the definition of $J$ in Theorem \ref{theorem3.4},
shows that 
\begin{equation}
J\left( w\right) =\frac{1}{2}\left\Vert 1_{\omega }^{\ast }e^{\left(
T_{3}-T_{2}\right) A}w\right\Vert _{\omega }^{2}+\left\langle z,e^{\left(
T_{3}-T_{1}\right) A}w\right\rangle +\varepsilon \left\Vert z\right\Vert
\left\Vert w\right\Vert =-\frac{1}{2}\left\Vert f\right\Vert _{\omega }^{2}%
\text{ .}  \tag{3.3.5}  \label{3.3.5}
\end{equation}%
Meanwhile, it follows from the above and the observation estimate in $(iii)$
in Theorem \ref{theorem3.8} that%
\begin{equation}
\begin{array}{ll}
J\left( w\right)  & \geq \frac{1}{2}\left\Vert 1_{\omega }^{\ast }e^{\left(
T_{3}-T_{2}\right) A}w\right\Vert _{\omega }^{2}+\varepsilon \left\Vert
z\right\Vert \left\Vert w\right\Vert -\left\Vert \,e^{\left(
T_{3}-T_{1}\right) A}w\right\Vert \left\Vert z\right\Vert  \\ 
& \geq \frac{1}{2}\left\Vert 1_{\omega }^{\ast }e^{\left( T_{3}-T_{2}\right)
A}w\right\Vert _{\omega }^{2}+\varepsilon \left\Vert z\right\Vert \left\Vert
w\right\Vert -\left( C\left\Vert 1_{\omega }^{\ast }e^{\left(
T_{3}-T_{2}\right) A}w\right\Vert _{\omega }+\varepsilon \left\Vert
w\right\Vert \right) \left\Vert z\right\Vert  \\ 
& \geq \frac{1}{2}\left\Vert f\right\Vert _{\omega }^{2}-C\left\Vert
f\right\Vert _{\omega }\left\Vert z\right\Vert \text{ .}%
\end{array}
\tag{3.3.6}  \label{3.3.6}
\end{equation}%
From (\ref{3.3.5}) and (\ref{3.3.6}), it follows that 
\begin{equation}
\left\Vert f\right\Vert _{\omega }\leq C\left\Vert z\right\Vert \text{ .} 
\tag{3.3.7}  \label{3.3.7}
\end{equation}%
On the other hand, since $f$ is the minimal norm control to $(\mathcal{P})$,
it holds that $\left\Vert y\left( T_{3}\right) \right\Vert \leq \varepsilon
\left\Vert z\right\Vert $. From this and (\ref{3.3.7}), we find that (\ref%
{3.3.2}) is true. Hence, $(ii)$ stands.

In summary, we complete the proof of Theorem \ref{theorem3.8}.

\bigskip

The rest of this subsection is devoted to the proof of Theorem \ref%
{theorem3.7}.

\bigskip

Proof of Theorem \ref{theorem3.7} .- Define 
\begin{equation*}
C^{\ast }:=\text{inf}\left\{ C>0;C\text{ satisfies (\ref{3.3.2})}\right\} 
\text{ .}
\end{equation*}%
\noindent \textit{Proof of} $\mathcal{N}\leq C^{\ast }$: It directly follows
from Theorem \ref{theorem3.8}.

\noindent \textit{Proof of} $C^{\ast }\leq \mathcal{N}$: It suffices to
prove that the property $(\mathcal{Q}_{\mathcal{N}})$ holds. Indeed, by
making use of the proof of \textquotedblleft $(ii)\Rightarrow \left(
iii\right) $\textquotedblright\ of Theorem \ref{theorem3.8}, we find that 
\begin{equation*}
C^{\ast }\leq \mathcal{N}\text{ .}
\end{equation*}%
Therefore, the remainder is that for each $z\in L^{2}\left( \Omega \right) $%
, there is a control $f\in L^{2}\left( \omega \right) $ satisfying that%
\begin{equation}
\text{max}\left\{ \frac{1}{\mathcal{N}}\left\Vert f\right\Vert _{\omega },%
\frac{1}{\varepsilon }\left\Vert y\left( T_{3}\right) \right\Vert \right\}
\leq \left\Vert z\right\Vert \text{ .}  \tag{3.3.8}  \label{3.3.8}
\end{equation}%
When $z=0$, we can easily get (\ref{3.3.8}) by taking $f=0$. So it suffices
to prove (\ref{3.3.8}) for an arbitrarily fixed $z\in L^{2}\left( \Omega
\right) \left\backslash \left\{ 0\right\} \right. $. To this end, we let $%
\widehat{z}=z/\left\Vert z\right\Vert $. Denote $y$ the solution of (\ref%
{3.2.1}) associated with the initial datum $\widehat{z}$. It follows from (%
\ref{3.3.1}) that 
\begin{equation*}
\text{inf}\left\{ \left\Vert f\right\Vert _{\omega };f\in L^{2}\left( \omega
\right) \text{ and }\left\Vert y\left( T_{3}\right) \right\Vert \leq
\varepsilon \left\Vert \widehat{z}\right\Vert \right\} \leq \mathcal{N}\text{
.}
\end{equation*}%
Because the infimum on the left hand side of the above inequality can be
reached, there is $\widehat{f}\in L^{2}\left( \omega \right) $ so that%
\begin{equation*}
\left\Vert \widehat{y}\left( T_{3}\right) \right\Vert \leq \delta \left\Vert 
\widehat{z}\right\Vert \quad \text{with }\left\Vert \widehat{f}\right\Vert
_{\omega }\leq \mathcal{N}\text{ .}
\end{equation*}%
where $\widehat{y}$ the solution of (\ref{3.2.1}) associated with the
initial datum $\widehat{z}$ and control $\widehat{f}$. From these, we see
that (\ref{3.3.8}) holds for $f=\left\Vert z\right\Vert \widehat{f}$. This
ends the proof.

\bigskip

\bigskip

\bigskip

\section{Inverse source problem}

\bigskip

This section concerns an inverse source problem: Suppose that we have a
solution $\varphi $ of $\varphi ^{\prime }-A\varphi =0$ with a priori bound
on the initial data in $L^{2}(\Omega )$. The question is how to recover
approximatively the initial data from the knowledge of the solution $\varphi 
$ in the future. This can be done as follows thanks to the impulse control.

\bigskip

\begin{thm}
\label{theorem4.1} Let $\omega _{1}$ be a non-empty open subset of $\Omega $
and $K\in \mathbb{N}$. Let $0\leq T_{1}<T_{2}<T_{3}$ and let $\varphi $ be a
solution of 
\begin{equation*}
\left\{ 
\begin{array}{ll}
\varphi ^{\prime }\left( t\right) -A\varphi \left( t\right) =0\text{ ,} & 
t\in \left( T_{1},T_{3}\right) \ \text{,} \\ 
\varphi \left( T_{1}\right) \in L^{2}\left( \Omega \right) \text{ .} & 
\end{array}%
\right.
\end{equation*}%
Then for any $\varepsilon >0$, there exists $\left\{ g_{j}\right\}
_{j=1,\cdot \cdot ,K}\in L^{2}\left( \omega _{1}\right) $ such that for any $%
j=1,\cdot \cdot ,K$, 
\begin{equation*}
\left\vert \langle \varphi \left( T_{1}\right) ,\xi _{j}\rangle +e^{\left(
T_{3}-T_{1}\right) \lambda _{j}}\langle g_{j},1_{\omega _{1}}^{\ast }\varphi
\left( T_{1}+T_{3}-T_{2}\right) \rangle _{\omega _{1}}\right\vert \leq
e^{\left( T_{3}-T_{1}\right) \lambda _{j}}\varepsilon \left\Vert \varphi
\left( T_{1}\right) \right\Vert
\end{equation*}%
and 
\begin{equation*}
\left\Vert g_{j}\right\Vert _{\omega _{1}}\leq \mathcal{C}_{\varepsilon
}\left( T_{3}-T_{2},T_{2}-T_{1}\right)
\end{equation*}%
where $\mathcal{C}_{\varepsilon }$ is given in Theorem \ref{theorem3.1}.
Further $g_{j}$ is the control function given in Corollary \ref{corollary3.2}
with $\omega =\omega _{1}$.
\end{thm}

\bigskip

Proof .-

\noindent \textit{Step 1:} We apply Corollary \ref{corollary3.2} with $%
\omega =\omega _{1}$ and get the existence of $\left( y_{j},g_{j}\right) $
such that 
\begin{equation*}
\left\{ 
\begin{array}{ll}
y_{j}^{\prime }\left( t\right) -Ay_{j}\left( t\right) =0\text{ ,} & t\in
\left( T_{1},T_{3}\right) \backslash \left\{ T_{2}\right\} \ \text{,} \\ 
y_{j}\left( T_{1}\right) =\xi _{j}\text{ ,} &  \\ 
y_{j}\left( T_{2}\right) =y_{j}\left( T_{2-}\right) +1_{\omega _{1}}g_{j}%
\text{ ,} & 
\end{array}%
\right.
\end{equation*}%
and $\left\Vert y_{j}\left( T_{3}\right) \right\Vert \leq \varepsilon $
where $g_{j}$ has the desired bound.

\noindent \textit{Step 2:} Write $\varphi \left( T_{1}\right) =\displaystyle%
\sum_{i=1,\cdot \cdot ,+\infty }a_{i}\xi _{i}$ with $a_{i}=\langle \varphi
\left( T_{1}\right) ,\xi _{i}\rangle $. Then we have that 
\begin{equation*}
\varphi \left( T_{3}\right) =\displaystyle\sum_{i=1,\cdot \cdot ,+\infty
}a_{i}e^{-\left( T_{3}-T_{1}\right) \lambda _{j}}\xi _{i}\text{ .}
\end{equation*}%
Hence, $\langle y_{j}(T_{1}),\varphi \left( T_{3}\right) \rangle =\langle
\xi _{j},\varphi \left( T_{3}\right) \rangle =a_{j}e^{-\left(
T_{3}-T_{1}\right) \lambda _{j}}$.

\noindent \textit{Step 3:} Multiply the equation solved by $y_{j}$ by the
solution $\varphi \left( T_{1}+T_{3}-t\right) $, with $t\in \left[
T_{1},T_{3}\right] $, to get 
\begin{equation*}
\langle y_{j}(T_{3}),\varphi \left( T_{1}\right) \rangle =\langle
y_{j}(T_{1}),\varphi \left( T_{3}\right) \rangle +\langle g_{j},1_{\omega
_{1}}^{\ast }\varphi \left( T_{1}+T_{3}-T_{2}\right) \rangle _{\omega _{1}}%
\text{ .}
\end{equation*}%
Therefore by step 2, it holds that 
\begin{equation*}
\left\vert a_{j}+e^{\left( T_{3}-T_{1}\right) \lambda _{j}}\langle
g_{j},1_{\omega _{1}}^{\ast }\varphi \left( T_{1}+T_{3}-T_{2}\right) \rangle
_{\omega _{1}}\right\vert =e^{\left( T_{3}-T_{1}\right) \lambda
_{j}}\left\vert \langle y_{j}(T_{3}),\varphi \left( T_{1}\right) \rangle
\right\vert \text{ .}
\end{equation*}%
This, along with the Cauchy-Schwarz inequality and step 1, leads to the
desired result.

\bigskip

\bigskip

\section{Main result}

\bigskip

This section presents the main result of this paper, as well as its proof.
We first recall that $\omega _{1}$ and $\omega _{2}$ are two arbitrarily
fixed open and non-empty subsets of $\Omega $. We next recall that $%
\{\lambda _{j}\}_{j=1}^{\infty }$ is the family of all eigenvalues of $-A$
so that (\ref{1.2}) holds and that $\{\xi _{j}\}_{j=1}^{\infty }$ is the
family of the corresponding normalized eigenfunctions. For each $\gamma >0$,
we define a natural number $K$ in the following manner: 
\begin{equation}
K:=\text{card}\left\{ j\in \mathbb{N},\text{ }\lambda _{j}<\gamma +\frac{%
\text{ln}2}{T}\right\} \text{ .}  \tag{5.1}  \label{5.1}
\end{equation}%
Next, we define 
\begin{equation}
\varepsilon :=\frac{1}{6\left( 1+K\right) }e^{-\gamma T}e^{-\left\Vert
V\right\Vert _{\infty }T}e^{-\lambda _{K}T/2}\text{ .}  \tag{5.2}
\label{5.2}
\end{equation}%
Denote by $\left\{ f_{j}\right\} _{j=1,\cdot \cdot ,K}\in L^{2}\left( \omega
_{2}\right) $ the minimal norm control functions obtained by applying
Corollary \ref{corollary3.2} with $T_{1}=\frac{T}{4}$, $T_{2}=T$, $T_{3}=%
\frac{5T}{4}$ and $\omega =\omega _{2}$. Denote by $\left\{ g_{j}\right\}
_{j=1,\cdot \cdot ,K}\in L^{2}\left( \omega _{1}\right) $ the minimal norm
control functions obtained by applying Corollary \ref{corollary3.2} with $%
T_{1}=\frac{T}{4}$, $T_{2}=\frac{T}{2}$, $T_{3}=\frac{3T}{4}$ and $\omega
=\omega _{1}$. Now for each $\gamma >0$, we define a linear bounded operator 
$\mathcal{F}$ from $L^{2}\left( \omega _{1}\right) $ into $L^{2}\left(
\omega _{2}\right) $ in the following manner: 
\begin{equation}
\mathcal{F}(p):=-\sum_{j=1,\cdot \cdot ,K}e^{\lambda _{j}T/2}\langle
g_{j},p\rangle _{\omega _{1}}f_{j}\left( x\right) \quad \text{for each }p\in
L^{2}(\omega _{1})\text{ .}  \tag{5.3}  \label{5.3}
\end{equation}%
The closed-loop equation under consideration reads: 
\begin{equation}
\left\{ 
\begin{array}{ll}
y^{\prime }\left( t\right) -Ay\left( t\right) =0\text{ ,} & \text{in }\left(
0,+\infty \right) \backslash \mathbb{N}T\ \text{,} \\ 
y\left( 0\right) \in L^{2}\left( \Omega \right) \text{ ,} &  \\ 
y\left( \left( n+1\right) T\right) =y\left( \left( n+1\right) T_{-}\right)
+1_{\omega _{2}}\mathcal{F}\left( 1_{\omega _{1}}^{\ast }y\left( \left( n+%
\frac{1}{2}\right) T\right) \right) \text{ ,} & \text{for }n\in \overline{%
\mathbb{N}}\text{ .}%
\end{array}%
\right.  \tag{5.4}  \label{5.4}
\end{equation}

\bigskip

The main result of this paper is the following theorem:

\bigskip

\begin{thm}
\label{theorem5.1} For each $\gamma >0$, let $\mathcal{F}$ be given by (\ref%
{5.3}). Then the following conclusions are true:

\begin{description}
\item[$\left( i\right) $] Each solution $y$ to the equation (\ref{5.4})
satisfies that 
\begin{equation*}
\left\Vert y\left( t\right) \right\Vert \leq e^{T\left( \gamma +\left\Vert
V\right\Vert _{\infty }\right) }\left( 1+\left\Vert \mathcal{F}\right\Vert _{%
\mathcal{L}(L^{2}(\omega _{1}),L^{2}(\omega _{2}))}\right) e^{-\gamma
t}\left\Vert y(0)\right\Vert \quad \text{for all}\; \;t\geq 0\text{ .}
\end{equation*}

\item[$\left( ii\right) $] The operator $\mathcal{F}$ satisfies the
estimate: 
\begin{equation*}
\left\Vert \mathcal{F}\right\Vert _{\mathcal{L}(L^{2}(\omega
_{1}),L^{2}(\omega _{2}))}\leq Ce^{C\gamma }\text{ ,}
\end{equation*}%
where $C$ is a positive constant independent of $\gamma $, depending on $%
\Omega $, $\omega _{1}$ ,$\omega _{2}$, $d$, $T$ and $\Vert V\Vert _{\infty
} $. Moreover, the manner how it depends on $T$, $d$, $\Vert V\Vert _{\infty
}$ is explicitly given.
\end{description}
\end{thm}

\bigskip

Proof .-

\noindent \textit{Step 1:}

Set $L_{n}=nT+\frac{T}{4}$. In order to have the conclusion $(i)$ in the
theorem, it suffices to prove that the solution $y$ of (\ref{5.4}) satisfies 
\begin{equation*}
\left\Vert y\left( L_{n+1}\right) \right\Vert \leq e^{-\gamma T}\left\Vert
y\left( L_{n}\right) \right\Vert
\end{equation*}%
for any $n\geq 0$. Indeed, thanks to the above inequality, we find that when 
$t\in \left[ L_{n},\left( n+1\right) T\right] $, 
\begin{equation*}
\begin{array}{ll}
\left\Vert y\left( t\right) \right\Vert & \leq e^{\left( t-L_{n}\right)
\left\Vert V\right\Vert _{\infty }}\left\Vert y\left( L_{n}\right)
\right\Vert \leq e^{\left( t-L_{n}\right) \left\Vert V\right\Vert _{\infty
}}e^{-n\gamma T}\left\Vert y(0)\right\Vert \\ 
& \leq e^{T\left\Vert V\right\Vert _{\infty }}e^{-n\gamma T}\left\Vert
y(0)\right\Vert \leq e^{T\left( \gamma +\left\Vert V\right\Vert _{\infty
}\right) }e^{-\gamma t}\left\Vert y(0)\right\Vert \text{ ,}%
\end{array}%
\end{equation*}%
and when $t\in \left[ \left( n+1\right) T,L_{n+1}\right] $, 
\begin{equation*}
\begin{array}{ll}
\left\Vert y\left( t\right) \right\Vert & \leq e^{\left( t-\left( n+1\right)
T\right) \left\Vert V\right\Vert _{\infty }}\left\Vert y\left( \left(
n+1\right) T\right) \right\Vert \\ 
& \leq e^{\left( t-\left( n+1\right) T\right) \left\Vert V\right\Vert
_{\infty }}\left\Vert y\left( \left( n+1\right) T\right) _{-}\right\Vert
+e^{\left( t-\left( n+1\right) T\right) \left\Vert V\right\Vert _{\infty
}}\left\Vert \mathcal{F}\left( 1_{\omega _{1}}^{\ast }y\left( \left( n+\frac{%
1}{2}\right) T\right) \right) \right\Vert _{\omega _{2}} \\ 
& \leq e^{\left( t-\left( n+1\right) T+3T/4\right) \left\Vert V\right\Vert
_{\infty }}\left\Vert y\left( L_{n}\right) \right\Vert +\left\Vert \mathcal{F%
}\right\Vert _{\mathcal{L}(L^{2}(\omega _{1}),L^{2}(\omega ))}e^{\left(
t-\left( n+1\right) T+T/4\right) \left\Vert V\right\Vert _{\infty
}}\left\Vert y\left( L_{n}\right) \right\Vert \\ 
& \leq e^{T\left( \gamma +\left\Vert V\right\Vert _{\infty }\right) }\left(
1+\left\Vert \mathcal{F}\right\Vert _{\mathcal{L}(L^{2}(\omega
_{1}),L^{2}(\omega _{2}))}\right) e^{-\gamma t}\left\Vert y(0)\right\Vert 
\text{ .}%
\end{array}%
\end{equation*}%
From these and the time translation invariance of the equation (\ref{5.4}),
we see that the conclusion $(i)$ in Theorem \ref{theorem5.1} is true for any 
$t\geq \frac{T}{4}$. But the case $t\leq \frac{T}{4}$ is trivial.

\noindent \textit{Step 2:}

Denote $y\left( L_{n}\right) =\displaystyle\sum_{j=1,\cdot \cdot ,+\infty
}a_{j}\xi _{j}$ and $a=\left( a_{j}\right) _{j=1,\cdot \cdot ,+\infty }$.
Then one deduces that $\langle y\left( L_{n}\right) ,\xi _{j}\rangle =a_{j}$
and $\left\Vert a\right\Vert _{\ell ^{2}}=\left\Vert y\left( L_{n}\right)
\right\Vert $.

For the rest of the proof, recall that $K$ and $\varepsilon $ are given by (%
\ref{5.1}) and (\ref{5.2}) respectively.

\noindent \textit{Step 3:}

Notice that the solution $y$ of (\ref{5.4}) evolves freely without a control
function between in $\left[ L_{n},nT+\frac{3T}{4}\right] $. Thus, we can
apply Theorem \ref{theorem4.1} with the choice $T_{1}=L_{n}$, $T_{2}=nT+%
\frac{T}{2}$, $T_{3}=nT+\frac{3T}{4}$ and $\varphi =y$ to get $\left\{
g_{j}\right\} _{j=1,\cdot \cdot ,K}\in L^{2}\left( \omega _{1}\right) $ such
that for any $j=1,\cdot \cdot ,K$, 
\begin{equation*}
\begin{array}{ll}
\left\vert a_{j}+e^{\lambda _{j}T/2}\langle g_{j},1_{\omega _{1}}^{\ast
}y\left( \left( n+\frac{1}{2}\right) T\right) \rangle _{\omega
_{1}}\right\vert & \leq e^{\lambda _{j}T/2}\varepsilon \left\Vert y\left(
L_{n}\right) \right\Vert \\ 
& \leq e^{\lambda _{K}T/2}\varepsilon \left\Vert y\left( L_{n}\right)
\right\Vert%
\end{array}%
\end{equation*}%
and 
\begin{equation*}
\left\Vert g_{j}\right\Vert _{\omega _{1}}\leq \mathcal{C}_{\varepsilon
}\left( T/4,T/4\right)
\end{equation*}%
where $\mathcal{C}_{\varepsilon }$ is given in Theorem \ref{theorem3.1}.
Further such $g_{j}$ is given in Corollary \ref{corollary3.2} with $\omega
=\omega _{1}$. By the time translation invariance of the equation (\ref{5.4}%
), $\left\{ g_{j}\right\} _{j=1,\cdot \cdot ,K}\in L^{2}\left( \omega
_{1}\right) $ is the control function obtained by applying Corollary \ref%
{corollary3.2} with $T_{1}=\frac{T}{4}$, $T_{2}=\frac{T}{2}$, $T_{3}=\frac{3T%
}{4}$ and $\omega =\omega _{1}$.

\noindent \textit{Step 4:}

Denote $b_{j}:=-e^{\lambda _{j}T/2}\langle g_{j},1_{\omega _{1}}^{\ast
}y\left( \left( n+\frac{1}{2}\right) T\right) \rangle _{\omega _{1}}$ for $%
j=1,\cdot \cdot ,K$, and $b=\left( b_{j}\right) _{j=1,\cdot \cdot ,K}$. Then
by step 3 and step 2, we have 
\begin{equation*}
\begin{array}{ll}
\left\Vert b\right\Vert _{\ell ^{2}} & \leq \left\Vert \left(
a_{j}+e^{\lambda _{j}T/2}\langle g_{j},1_{\omega _{1}}^{\ast }y\left( \left(
n+\frac{1}{2}\right) T\right) \rangle _{\omega _{1}}\right) _{j=1,\cdot
\cdot ,K}\right\Vert _{\ell ^{2}}+\left\Vert a\right\Vert _{\ell ^{2}} \\ 
& \leq \left( \sqrt{K}e^{\lambda _{K}T/2}\varepsilon +1\right) \left\Vert
y\left( L_{n}\right) \right\Vert \text{ .}%
\end{array}%
\end{equation*}

\noindent \textit{Step 5:}

We apply Theorem \ref{theorem3.3} with $T_{1}=L_{n}$, $T_{2}=\left(
n+1\right) T$, $T_{3}=L_{n+1}$ and the above choice of $b$. Then there is a
solution $\left( \widetilde{y},\widetilde{f}\right) $ such that 
\begin{equation*}
\left\{ 
\begin{array}{ll}
\widetilde{y}^{\prime }\left( t\right) -A\widetilde{y}\left( t\right) =0%
\text{ ,} & t\in \left( L_{n},L_{n+1}\right) \backslash \left\{ \left(
n+1\right) T\right\} \ \text{,} \\ 
\widetilde{y}\left( L_{n}\right) =\displaystyle\sum_{j=1,\cdot \cdot
,K}b_{j}\xi _{j}\text{ ,} &  \\ 
\widetilde{y}\left( \left( n+1\right) T\right) =\widetilde{y}\left( \left(
n+1\right) T_{-}\right) +1_{\omega _{2}}\widetilde{f}\text{ ,} & 
\end{array}%
\right.
\end{equation*}%
and 
\begin{equation*}
\left\{ 
\begin{array}{ll}
\left\Vert \widetilde{y}\left( L_{n+1}\right) \right\Vert \leq \varepsilon 
\sqrt{K}\left\Vert b\right\Vert _{\ell ^{2}}\text{ ,} &  \\ 
\widetilde{f}=\displaystyle\sum_{j=1,\cdot \cdot ,K}b_{j}f_{j}:=-%
\displaystyle\sum_{j=1,\cdot \cdot ,K}e^{\lambda _{j}T/2}\langle
g_{j},1_{\omega _{1}}^{\ast }y((n+\tfrac{1}{2})T)\rangle _{\omega _{1}}f_{j}%
\text{ ,} & 
\end{array}%
\right.
\end{equation*}%
where $f_{j}$ is given by Corollary \ref{corollary3.2} with $\omega =\omega
_{2}$ and satisfies%
\begin{equation*}
\left\Vert f_{j}\right\Vert _{\omega _{2}}\leq \mathcal{C}_{\varepsilon
}\left( T/4,3T/4\right)
\end{equation*}%
with $\mathcal{C}_{\varepsilon }$ given in Theorem \ref{theorem3.1}. By the
time translation invariance of the equation (\ref{5.4}), $\left\{
f_{j}\right\} _{j=1,\cdot \cdot ,K}\in L^{2}\left( \omega _{2}\right) $ is
the control function obtained by applying Corollary \ref{corollary3.2} with $%
T_{1}=\frac{T}{4}$, $T_{2}=T$, $T_{3}=\frac{5T}{4}$ and $\omega =\omega _{2}$%
.

\noindent \textit{Step 6:}

One can check that $\widetilde{f}=\mathcal{F}\left( 1_{\omega _{1}}^{\ast
}y\left( \left( n+\frac{1}{2}\right) T\right) \right) $ and $y=\widetilde{y}+%
\widehat{y}+\overline{y}$ where $\widetilde{y}$ is given in step 5, $%
\widehat{y}$ solves 
\begin{equation*}
\left\{ 
\begin{array}{ll}
\widehat{y}^{\prime }\left( t\right) -A\widehat{y}\left( t\right) =0\text{ ,}
& t\in \left( L_{n},L_{n+1}\right) \backslash \left\{ \left( n+1\right)
T\right\} \ \text{,} \\ 
\widehat{y}\left( L_{n}\right) =\displaystyle\sum_{j=1,\cdot \cdot ,K}\left(
a_{j}-b_{j}\right) \xi _{j}\text{ ,} &  \\ 
\widehat{y}\left( \left( n+1\right) T\right) =\widehat{y}\left( \left(
n+1\right) T_{-}\right) \text{ ,} & 
\end{array}%
\right.
\end{equation*}%
and $\overline{y}$ satisfies 
\begin{equation*}
\left\{ 
\begin{array}{ll}
\overline{y}^{\prime }\left( t\right) -A\overline{y}\left( t\right) =0\text{
,} & t\in \left( L_{n},L_{n+1}\right) \backslash \left\{ \left( n+1\right)
T\right\} \ \text{,} \\ 
\overline{y}\left( L_{n}\right) =\displaystyle\sum_{j>K}a_{j}\xi _{j}\text{ ,%
} &  \\ 
\overline{y}\left( \left( n+1\right) T\right) =\overline{y}\left( \left(
n+1\right) T_{-}\right) \text{ .} & 
\end{array}%
\right.
\end{equation*}

\noindent \textit{Step 7:} We estimate $\left\Vert y\left( L_{n+1}\right)
\right\Vert =\left\Vert \widetilde{y}\left( L_{n+1}\right) +\widehat{y}%
\left( L_{n+1}\right) +\overline{y}\left( L_{n+1}\right) \right\Vert $ as
follows.

First, by step 5 and step 4, 
\begin{equation*}
\left\Vert \widetilde{y}\left( L_{n+1}\right) \right\Vert \leq \varepsilon 
\sqrt{K}\left\Vert b\right\Vert _{\ell ^{2}}\leq \varepsilon \sqrt{K}\left( 
\sqrt{K}e^{\lambda _{K}T/2}\varepsilon +1\right) \left\Vert y\left(
L_{n}\right) \right\Vert \text{ .}
\end{equation*}%
Second, by step 3%
\begin{equation*}
\left\Vert \widehat{y}\left( L_{n+1}\right) \right\Vert \leq e^{\left\Vert
V\right\Vert _{\infty }T}\left\Vert \widehat{y}\left( L_{n}\right)
\right\Vert \leq e^{\left\Vert V\right\Vert _{\infty }T}\sqrt{K}e^{\lambda
_{K}T/2}\varepsilon \left\Vert y\left( L_{n}\right) \right\Vert \text{ .}
\end{equation*}%
Third, 
\begin{equation*}
\left\Vert \overline{y}\left( L_{n+1}\right) \right\Vert =\left(
\sum_{j>K}\left\vert a_{j}e^{-\lambda _{j}\left( L_{n+1}-L_{n}\right)
}\right\vert ^{2}\right) ^{1/2}\leq e^{-\lambda _{K+1}T}\left\Vert y\left(
L_{n}\right) \right\Vert \text{ .}
\end{equation*}%
Gathering all the previous estimates, one concludes that 
\begin{equation*}
\left\Vert y\left( L_{n+1}\right) \right\Vert \leq \left( e^{-\lambda
_{K+1}T}+3e^{\left\Vert V\right\Vert _{\infty }T}e^{\lambda _{K}T/2}\left(
1+K\right) \varepsilon \right) \left\Vert y\left( L_{n}\right) \right\Vert
\end{equation*}%
with $\varepsilon \in \left( 0,1\right) $. Finally, the choice of $K$ (see (%
\ref{5.1})) gives $e^{-\lambda _{K+1}T}\leq \frac{1}{2}e^{-\gamma T}$, and
the choice of $\varepsilon \in \left( 0,1\right) $ (see (\ref{5.2})) gives $%
3e^{\left\Vert V\right\Vert _{\infty }T}e^{\lambda _{K}T/2}\left( 1+K\right)
\varepsilon =\frac{1}{2}e^{-\gamma T}$, which implies the desired estimate
for $\left\Vert y\left( L_{n+1}\right) \right\Vert $, that is $\left\Vert
y\left( L_{n+1}\right) \right\Vert \leq e^{-\gamma T}\left\Vert y\left(
L_{n}\right) \right\Vert $.

\noindent \textit{Step 8:} We treat the boundedness of $\mathcal{F}$ as
follows:%
\begin{equation*}
\begin{array}{ll}
\left\Vert \mathcal{F}\left( w\right) \right\Vert _{\omega _{2}}^{2} & =%
\displaystyle\int_{\omega _{2}}\left\vert \displaystyle\sum_{j=1,\cdot \cdot
,K}e^{\lambda _{j}T/2}\langle g_{j},w\rangle _{\omega _{1}}f_{j}\left(
x\right) \right\vert ^{2}dx \\ 
& \leq \displaystyle\sum_{j=1,\cdot \cdot ,K}\left\vert e^{\lambda
_{j}T/2}\langle g_{j},w\rangle _{\omega _{1}}\right\vert ^{2}\displaystyle%
\sum_{j=1,\cdot \cdot ,K}\displaystyle\int_{\omega _{2}}\left\vert
f_{j}\left( x\right) \right\vert ^{2}dx \\ 
& \leq \left\Vert w\right\Vert _{\omega _{1}}^{2}e^{\lambda _{K}T}%
\displaystyle\sum_{j=1,\cdot \cdot ,K}\left\Vert g_{j}\right\Vert _{\omega
_{1}}^{2}\displaystyle\sum_{j=1,\cdot \cdot ,K}\left\Vert f_{j}\right\Vert
_{\omega _{2}}^{2} \\ 
& \leq \left\Vert w\right\Vert _{\omega _{1}}^{2}e^{\lambda _{K}T}K^{2}\left[
\mathcal{C}_{\varepsilon }\left( T/4,T/4\right) \mathcal{C}_{\varepsilon
}\left( T/4,3T/4\right) \right] ^{2}\text{ }%
\end{array}%
\end{equation*}%
which implies 
\begin{equation*}
\left\Vert \mathcal{F}\right\Vert _{\mathcal{L}(L^{2}(\omega
_{1}),L^{2}(\omega _{2}))}\leq e^{\lambda _{K}T/2}K\mathcal{C}_{\varepsilon
}\left( T/4,T/4\right) \mathcal{C}_{\varepsilon }\left( T/4,3T/4\right) 
\text{ .}
\end{equation*}

\noindent \textit{Step 9:} We estimate $1/\varepsilon $ and $\mathcal{C}%
_{\varepsilon }$:

By the Weyl's asymptotic law for the Dirichlet eigenvalues $\alpha _{j}$,
there is a constant $\overline{C}>0$ (depending only on $\Omega $ and $d$)
such that for any $\mu >0$,%
\begin{equation*}
\text{card}\left\{ j\in \mathbb{N},\text{ }\alpha _{j}<\mu \right\} \leq 
\overline{C}\left( 1+\mu ^{d/2}\right) \text{ .}
\end{equation*}%
By the min-max formula, one has 
\begin{equation*}
-\left\Vert V\right\Vert _{\infty }+\alpha _{j}\leq \lambda _{j}\leq \alpha
_{j}+\left\Vert V\right\Vert _{\infty }\text{ . }
\end{equation*}%
One deduces that there is a constant $\overline{C}>0$ (depending only on $%
\Omega $ and $d$) such that for any $\gamma >0$ 
\begin{equation*}
K:=\text{card}\left\{ j\in \mathbb{N},\text{ }\lambda _{j}<\gamma +\frac{%
\text{ln}2}{T}\right\} \leq \overline{C}\left( 1+\left\Vert V\right\Vert
_{\infty }^{d/2}+\left( \frac{\text{ln}2}{T}\right) ^{d/2}+\gamma
^{d/2}\right) \text{ .}
\end{equation*}%
Further, for some constant $\overline{C}>0$ (depending only on $\Omega $ and 
$d$), we have 
\begin{equation*}
\frac{1}{\varepsilon }:=6e^{\gamma T}e^{\left\Vert V\right\Vert _{\infty
}T}e^{\lambda _{K}T/2}\left( 1+K\right) \leq \overline{C}\left( \gamma
^{d/2}+\left\Vert V\right\Vert _{\infty }^{d/2}+\frac{1}{T^{d/2}}\right)
e^{\left\Vert V\right\Vert _{\infty }T}e^{2\gamma T}\text{ .}
\end{equation*}%
We finish the proof by gathering the previous estimates with the definition
of $\mathcal{C}_{\varepsilon }$, that is $\mathcal{C}_{\varepsilon }\left(
t,s\right) :=e^{4s\left\Vert V\right\Vert _{\infty }}e^{c\left( 1+\frac{1}{t}%
+t\left\Vert V\right\Vert _{\infty }+\left\Vert V\right\Vert _{\infty
}^{2/3}\right) }$exp$\left( \sqrt{\frac{c}{t}\text{ln}^{+}\frac{1}{%
\varepsilon }}\right) $.

Hence, we complete the proof of Theorem~\ref{theorem5.1}.

\bigskip

\bigskip

\bigskip

\bigskip


\bigskip

\bigskip

\bigskip

\bigskip

\begin{thebibliography}{AEWZ}
\bibitem[AEWZ]{AEWZ} J. Apraiz, L. Escauriaza, G. Wang and C. Zhang,
Observability inequalities and measurable sets. J. Eur. Math. Soc. (JEMS) 16
(2014), no. 11, 2433--2475.

\bibitem[Ba]{Ba} V. Barbu, Stabilization of Navier-Stokes flows,
Communications and Control Engineering Series, Springer, London, 2011.

\bibitem[BC]{BC} M. Bardi and I. Capuzzo-Dolcetta, Optimal control and
viscosity solutions of Hamilton-Jacobi-Bellman equations, With appendices by
Maurizio Falcone and Pierpaolo Soravia, Systems \& Control: Foundations \&
Applications, Birkh\"{a}user Boston, Inc., Boston, MA, 1997.

\bibitem[Be]{Be} A. Bensoussan, Dynamic programming and inventory control,
Studies in Probability, Optimization and Statistics, 3, IOS Press,
Amsterdam, 2011.

\bibitem[BL1]{BL1} A. Bensoussan and J.-L. Lions, Nouvelles m\'{e}thodes en
contr\^{o}le impulsionnel (French), Appl. Math. and Optim., 1 (1974/75),
289--312.

\bibitem[BL2]{BL2} A. Bensoussan and J.-L. Lions, Contr\^{o}le impulsionnel
et in\'{e}quations quasi variationnelles (French), [Impulse control and
quasivariational inequalities], M\'{e}thodes Math\'{e}matiques de
l'Informatique [Mathematical Methods of Information Science], 11,
Gauthier-Villars, Paris, 1982.

\bibitem[Br]{Br} R. Brockett, A stabilization problem, Open problems in
mathematical systems and control theory, 75--78, Commun. Control Engrg.
Ser., Springer, London, 1999.

\bibitem[CCo]{CCo} E. Cerpa and J.-M. Coron, Rapid stabilization for a
Korteweg-de Vries equation from the left Dirichlet boundary condition, IEEE
Trans. Automat. Control, 58 (2013), 1688--1695.

\bibitem[CCr]{CCr} E. Cerpa and E. Cr\'{e}peau, Rapid exponential
stabilization for a linear Korteweg-de Vries equation, Discrete Contin. Dyn.
Syst. Ser. B, 11 (2009), 655--668.

\bibitem[C]{C} J.-M. Coron, Control and nonlinearity, Mathematical Surveys
and Monographs, 136, American Mathematical Society, Providence, RI, 2007.

\bibitem[Co]{Co} J.-M. Coron, On the stabilization of controllable and
observable systems by an output feedback law, Math. Control Signals Systems,
7 (1994), 187--216.

\bibitem[CL1]{CL1} J.-M. Coron and Q. L\"{u}, Fredholm transform and local
rapid stabilization for a Kuramoto-Sivashinsky equation, J. Differential
Equations, 259 (2015), 3683--3729.

\bibitem[CL2]{CL2} J.-M. Coron and Q. L\"{u}, Local rapid stabilization for
a Korteweg-de Vries equation with a Neumann boundary control on the right,
J. Math. Pures Appl., (9) 102 (2014), 1080--1120.

\bibitem[Cu]{Cu} R. F. Curtain, Stabilization of parabolic systems with
point observation and boundary control via integral dynamic output feedback
of a finite-dimensional compensator, Analysis and optimization of systems
(Versailles, 1982), 761--775, Lecture Notes in Control and Inform. Sci., 44,
Springer, Berlin, 1982.

\bibitem[DS]{DS} V. A. Dykhta and O. N. Samsonyuk, Optimal'noe impul'snoe
upravlenie s prilozheniyami (Russian), [Optimal impulse control with
applications], Fizmatlit "Nauka", Moscow, 2000.

\bibitem[EFV]{EFV} L. Escauriaza, F. J. Fern\'{a}ndez and S. Vessella,
Doubling properties of caloric functions, Appl. Anal., 85 (2006), 205--223.

\bibitem[FPZ]{FPZ} C. Fabre, J.-P. Puel and E. Zuazua, Approximate
controllability of the semilinear heat equation, Proc. Roy. Soc. Edinburgh
Sect. A, 125 (1995), 31--61.

\bibitem[I]{I} K. Ichikawa, Output feedback stabilization, Internat. J.
Control, 16 (1972), 513--522.

\bibitem[JL]{JL} D. Jerison and G. Lebeau, Nodal sets of sums of
eigenfunctions, Harmonic analysis and partial differential equations
(Chicago, IL, 1996), Univ. Chicago Press, 1999, pp. 223--239.

\bibitem[K]{K} V. Komornik, Rapid boundary stabilization of linear
distributed systems, SIAM J. Control Optim., 35 (1997), 1591--1613.

\bibitem[LR]{LR} G. Lebeau and L. Robbiano, Contr\^{o}le exact de l'\'{e}%
quation de la chaleur, Comm. Partial Differential Equations, 20 (1995),
335--356.

\bibitem[LZ]{LZ} G. Lebeau and E. Zuazua, Null-controllability of a system
of linear thermoelasticity. Arch. Rational Mech. Anal. 141 (1998), no. 4,
297--329.

\bibitem[LL]{LL} J. Le Rousseau and G. Lebeau, On Carleman estimates for
elliptic and parabolic operator, Applications to unique continuation and
control of parabolic equations, ESAIM: COCV, 18 (2012), 712--747.

\bibitem[LY]{LY} X. Li and J. Yong, Optimal control theory for
infinite-dimensional systems, Systems \& Control: Foundations \&
Applications, Birkh\"{a}user Boston, Inc., Boston, MA, 1995.

\bibitem[L]{L} F. H. Lin, A uniqueness theorem for the parabolic equation,
Comm. Pure Appl. Math., 43 (1990), 127--136.

\bibitem[Li]{Li} J.-L. Lions, Exact controllability, stabilization and
perturbations for distributed systems. SIAM Rev. 30 (1988), no. 1, 1--68.

\bibitem[LM]{LM} S. I. Lyashko and A. A. Man'kovski\v{i}, Gradientnye metody
v zadachakh optimal'nogo impul'snogo upravleniya dlya sistem s
raspredelennymi parametrami (Russian), [Gradient methods in problems of
optimal impulse control for systems with distributed parameters],
[Preprint], 83-9, Akad. Nauk Ukrain. SSR, Inst. Kibernet., Kiev, 1983.

\bibitem[MR]{MR} B. M. Miller and E. Y. Rubinovich, Impulsive control in
continuous and discrete-continuous systems, Kluwer Academic/Plenum
Publishers, New York, 2003.

\bibitem[NS]{NS} D. Ne\u{s}i\'{c} and E. D. Sontag, Input-to-state
stabilization of linear systems with positive outputs, Systems Control
Lett., 35 (1998), 245--255.

\bibitem[OS]{OS} B. \O ksendal and A. Sulem, Applied stochastic control of
jump diffusions, Universitext, Springer-Verlag, Berlin, 2005.

\bibitem[PW]{PW} K. D. Phung and G. Wang, Quantitative unique continuation
for the semilinear heat equation in a convex domain, J. Funct. Anal., 259
(2010), 1230--1247.

\bibitem[PW1]{PW1} K. D. Phung and G. Wang, An observability estimate for
parabolic equations from a measurable set in time and its applications, J.
Eur. Math. Soc. (JEMS), 15 (2013), 681--703.

\bibitem[PWZ]{PWZ} K. D. Phung, L. Wang and C. Zhang, Bang-bang property for
time optimal control of semilinear heat equation, Ann. Inst. H. Poincar\'{e}
Anal. Non Lin\'{e}aire, 31 (2014), 477--499.

\bibitem[R]{R} S. F. Richard, Optimal impulse control of a diffusion process
with both fixed and proportional costs of control, SIAM J. Control Optim.,
15 (1977), 79--91.

\bibitem[TWZ]{TWZ} E. Tr\'{e}lat, L. Wang and Y. Zhang, Impulse and
sampled-data optimal control of heat equations, and error estimates, SIAM J.
Control Optim., 57 (2016), 2787--2819.

\bibitem[U]{U} J. Urquiza, Rapid exponential feedback stabilization with
unbounded control operators, SIAM J. Control Optim., 43 (2005), 2233--2244.

\bibitem[V]{V} A. Vest, Rapid stabilization in a semigroup framework, SIAM
J. Control Optim., 51 (2013), 4169--4188.

\bibitem[YY]{YY} V. A. Yakubovich and E. D. Yakubovich, Optimal
stabilization of a control system with constraints on the output variable
(Russian), Avtomat. i Telemekh, 1993, 79--88; translation in Automat. Remote
Control, 54 (1993), part 1, 1368--1376 (1994).

\bibitem[Z]{Z} J. Zabczyk, Mathematical control theory: an introduction,
Systems \& Control: Foundations \& Applications, Birkh\"{a}user Boston,
Inc., Boston, MA, 1992.
\end{thebibliography}
\end{document}